\documentclass[10pt]{amsart}

\usepackage{amsmath}
\usepackage{amsthm}
\usepackage[usenames,dvipsnames,svgnames,table]{xcolor}
\usepackage{amssymb}
\begingroup
    \makeatletter
    \@for\theoremstyle:=definition,remark,plain\do{%
        \expandafter\g@addto@macro\csname th@\theoremstyle\endcsname{%
            \addtolength\thm@preskip\parskip
            }%
        }
\endgroup

\usepackage[margin=10pt,font=small,labelfont=bf]{caption}

\usepackage[pdfauthor=Lyons\ and\ White,bookmarksnumbered,hyperfigures,pdftex,pagebackref,colorlinks=true,linkcolor=BrickRed,urlcolor=BrickRed,citecolor=OliveGreen,pdfstartview=FitH]{hyperref}
\usepackage{microtype}
\usepackage{mathtools}
\mathtoolsset{centercolon}
\numberwithin{equation}{section}
\numberwithin{figure}{section}
\usepackage{GWmacros}
\usepackage[final]{showlabels}

\theoremstyle{plain}

\newcommand{\bp}{o}

\renewcommand{\defn}[1]{\textbf{\textit{#1}}}

\newcommand\I[1]{\mathbf1_{#1}}

\newcommand\ro{r_{\mathrm o}}
\newcommand\st{\,;\ }
\newcommand\rt{\mathbf{r}}   
\newcommand\dist{\mathrm{dist}}  
\newcommand\bd{\partial}  
\newcommand\wrd{\mathbf w}  
\newcommand\Exp{\mathrm{Exp}}  
\newcommand\lep{\preccurlyeq}

\newcommand\cf{C}   
\newcommand\cfs{\mathbf C}   
\newcommand\bits{\mathbf b}   
\newcommand\cfm{\xi}   
\newcommand\seqs{\mathbf{s}}   
\newcommand\soltn{\zeta}  

\makeatletter
\@namedef{subjclassname@2020}{%
  \textup{2020} Mathematics Subject Classification}
\makeatother

\let\P\relax
\DeclareMathOperator{\P}{\mathbf{P}\mathopen{}}
\DeclareMathOperator{\E}{\mathbf{E}\mathopen{}}

\def\Psub_#1{\P_{\! #1}}

\def\Psubbig_#1#2{\Psub_{#1}\mkern-1.5mu\bigl[#2\bigr]}

\def\Psubbigg_#1#2{\Psub_{#1}\mkern-1.5mu\biggl[#2\biggr]}

\def\Ebig#1{\E\mkern-1.5mu\bigl[#1\bigr]}
\def\Esubbig_#1#2{\E_{#1}\mkern-1.5mu\bigl[#2\bigr]}
\def\EsubBig_#1#2{\E_{#1}\mkern-1.5mu\Bigl[#2\Bigr]}
\def\Esupbig^#1#2{\E^{#1}\mkern-1.5mu\bigl[#2\bigr]}
\def\EsupBig^#1#2{\E^{#1}\mkern-1.5mu\Bigl[#2\Bigr]}

\def\Esubbigg_#1#2{\E_{#1}\mkern-1.5mu\biggl[#2\biggr]}

\def\thmenv#1#2#3{\begin{#1} \label{#1:#2} #3 \end{#1}}
\def\richthmenv#1#2#3#4{\begin{#1}[#3] \label{#1:#2} #4 \end{#1}}

\def\procl#1.#2 #3\endprocl{%
       \ifx#1t\thmenv{Theorem}{#2}{#3}\fi
       \ifx#1l\thmenv{Lemma}{#2}{#3}\fi
       \ifx#1p\thmenv{Proposition}{#2}{#3}\fi
       \ifx#1c\thmenv{Corollary}{#2}{#3}\fi
       \ifx#1d\thmenv{Definition}{#2}{#3}\fi
       \ifx#1g\thmenv{Conjecture}{#2}{#3}\fi
       \ifx#1q\thmenv{Question}{#2}{#3}\fi
       \ifx#1r\thmenv{Remark}{#2}{#3}\fi
    }%

\def\rprocl#1.#2 #3 #4\endprocl{%
       \ifx#1t\richthmenv{Theorem}{#2}{#3}{#4}\fi
       \ifx#1l\richthmenv{Lemma}{#2}{#3}{#4}\fi
       \ifx#1p\richthmenv{Proposition}{#2}{#3}{#4}\fi
       \ifx#1c\richthmenv{Corollary}{#2}{#3}{#4}\fi
       \ifx#1d\richthmenv{Definition}{#2}{#3}{#4}\fi
       \ifx#1g\richthmenv{Conjecture}{#2}{#3}{#4}\fi
       \ifx#1q\richthmenv{Question}{#2}{#3}{#4}\fi
       \ifx#1r\richthmenv{Remark}{#2}{#3}{#4}\fi
    }%

\def\rref#1.#2/{%
      \ifx #1sSection~\ref{s.#2}\fi
      \ifx #1SSubsection~\ref{S.#2}\fi
      \ifx #1tTheorem~\ref{Theorem:#2}\fi  
      \ifx #1lLemma~\ref{Lemma:#2}\fi 
      \ifx #1cCorollary~\ref{Corollary:#2}\fi 
      \ifx #1pProposition~\ref{Proposition:#2}\fi 
      \ifx #1dDefinition~\ref{Definition:#2}\fi
      \ifx #1gConjecture~\ref{Conjecture:#2}\fi 
      \ifx #1qQuestion~\ref{Question:#2}\fi 
      \ifx #1rRemark~\ref{Remark:#2}\fi 
      \ifx #1aAppendix~\ref{a.#2}\fi 
      \ifx #1fFigure~\ref{f.#2}\fi
      \ifx #1e(\ref{e.#2})\fi
      \ifx #1b\cite{#2}\fi
      \ifx #1B\cite{#2}\fi
        }

\def\briefref#1.#2/{%
      \ifx #1s\ref{s.#2}\fi
      \ifx #1S\ref{S.#2}\fi
      \ifx #1t\ref{Theorem:#2}\fi  
      \ifx #1l\ref{Lemma:#2}\fi 
      \ifx #1c\ref{Corollary:#2}\fi 
      \ifx #1p\ref{Proposition:#2}\fi 
      \ifx #1d\ref{Definition:#2}\fi
      \ifx #1g\ref{Conjecture:#2}\fi 
      \ifx #1q\ref{Question:#2}\fi 
      \ifx #1r\ref{Remark:#2}\fi 
      \ifx #1x\ref{Example:#2}\fi
      \ifx #1a\ref{a.#2}\fi 
      \ifx #1f\ref{f.#2}\fi
      \ifx #1e(\ref{e.#2})\fi
      \ifx #1b\cite{#2}\fi
      \ifx #1B\cite{#2}\fi
        }

\def\rlabel #1 #2{\begin{equation} \label{#1} #2 \end{equation}}

\def\rproof{\begin{proof}}

\def\Qed{\end{proof}}

\def\eqaln#1{\begin{align*} #1 \end{align*}}

\title{Monotonicity for continuous-time random walks}
\date{29 Sep.\ 2022}
\author{Russell Lyons and Graham White}
\address{Department of Mathematics, 831 E. 3rd St.,
Indiana University, Bloomington, IN 47405-7106} 
\email{\href{mailto:rdlyons@indiana.edu}{rdlyons@indiana.edu}}
\thanks{%
The work of R.L.\ was partially supported by the National
Science Foundation under grants DMS-1612363 and DMS-1954086 and also by the Simons Foundation.}

\email{\href{mailto:grahamwhite@alumni.stanford.edu}{grahamwhite@alumni.stanford.edu}}

\begin{document}

\begin{abstract}
Consider continuous-time random walks on Cayley graphs where the rate assigned to each edge depends only on the corresponding generator. We show that the limiting speed is monotone increasing in the rates for infinite Cayley graphs that arise from Coxeter systems, but not for all Cayley graphs. On finite Cayley graphs, we show that the distance --- in various senses --- to stationarity is monotone decreasing in the rates for Coxeter systems and for abelian groups, but not for all Cayley graphs. We also find several examples of surprising behaviour in the dependence of the distance to stationarity on the rates. This includes a counterexample to a conjecture on entropy of Benjamini, Lyons, and Schramm. We also show that the expected distance at any fixed time for random walks on $\Z^+$ is monotone increasing in the rates for arbitrary rate functions, which is not true on all of $\Z$. Various intermediate results are also of interest.
\end{abstract}



\maketitle

\section{Introduction}

We are interested in questions arising from the following scenario. Consider a connected, undirected graph $G$ with each edge $e$ labelled by a nonnegative real number $r(e)$. A random walk is defined on this graph by associating a Poisson process (``clock'') with rate $r(e)$ to each edge $e$. When this clock rings, if the walker is at either neighbouring vertex, then it moves along the edge $e$ to the other vertex; otherwise it does not move. 
We choose the convention that the walk is left-continuous in time, which will be more convenient to describe our couplings. 
Note that if $G$ is finite, then the stationary distribution of this walk is uniform.

Monotonicity in time or in the rate function $\rt := r(\cdot)$ of the behaviour of such a random walk is complicated and not always intuitive, primarily because increasing $r(e)$ not only makes the walk more likely to cross $e$ when at each of its endpoints, but also because it moves the walk sooner in both directions of traversing $e$. Let $p_t(x, y)$ be the transition probability from $x$ to $y$ at time $t$. Note that changing the rate function by a constant factor $\alpha$ is equivalent to changing the time $t$ to $\alpha \cdot t$. When $x = y$ and $t$ is fixed, there are examples where $p_t(x, x)$ is not monotone in the rate function (one such is given here in Example \ref{ex:Zdist}). Nevertheless, if $G$ is finite, then the average of $p_t(x, x)$ over all vertices $x$ of $G$ \emph{is} monotone, as shown by Benjamini and Schramm \cite[Theorem 3.1]{HeicklenHoffman}. In particular, if the graph with rates is vertex-transitive, then $p_t(x, x)$ is monotone decreasing in the rates. This result was extended in various ways by \cite{FontesMathieu,Lyons:treedom,AL:urn,PittSC,Lyons:avgrtn}, but all such extensions concerned only return probabilities. 

In a different direction, Karlin and McGregor \cite[equation (45)]{KarlinMcG} considered con\-tin\-uous-time random walks on $\Z^+ := \{0, 1, 2, \dots\}$ with arbitrary, symmetric rate functions.
It follows from their result
that for each $n \ge 1$, the time to reach $n$ from $0$ stochastically decreases in $\rt$. Indeed, they gave an explicit representation of that time
as the sum of $n$ independent, exponential random variables whose rates are the nonzero eigenvalues of the $(n+1) \times (n+1)$ tridiagonal Laplacian
matrix $\Delta$, where $\Delta(i, i) := r(i, i+1) + r(i, i-1)$, $\>\Delta(i, i+1) := -r(i, i+1)$, $\>\Delta(i, i-1) := -r(i, i-1)$, and all other entries
are $0$. Because $\Delta$ itself is monotone increasing in $\rt$ in the Loewner order, its eigenvalues are also increasing in the rate function.
(This monotonicity is also what is behind all the results mentioned above on the return probabilities. We note that \cite[equation (45)]{KarlinMcG} holds for more general birth and death chains; \cite{KarlinMcG:diff} gives a representation of transition probabilities via orthogonal polynomials.)
See also \cite{Fill:1st,DiaconisMiclo,Fill:general,Miclo} for other proofs and extensions of the representation \cite[equation (45)]{KarlinMcG}.

We are motivated by the following observation \cite[Exercise 13.24]{LP:book}: If the 3-regular tree is regarded as the standard Cayley graph of the free product of 3 copies of $\Z/2\Z$, and the rates on all edges are 1 except for those corresponding to a fixed generator, where the rates are $\rho$, then the limiting rate of escape (i.e., the graph distance from the starting point divided by the time as time tends to infinity) is 
\[
\bigl(3 \rho(\rho+1) + (1-\rho)\sqrt{16 \rho + 9 \rho^2}\, \bigr)\big/\bigl(2 (2 + \rho)\bigr) \text{ a.s.}
\]
What is interesting about this formula is that it is monotone increasing in $\rho$, but that this fact is not obvious. Is there a more intuitive explanation for such monotonicity? One might attempt to answer this via a coupling of two such random walks, where one is always at distance at least as great as the other, but no such coupling can be Markovian. In fact, even non-Markovian couplings cannot always have this property: see Proposition \ref{prop:cantforever}. 
Other questions of interest include these: For rates that depend only on the generators for the free product of more copies of $\Z/2\Z$, does monotonicity still hold? What about other Cayley graphs? What about other aspects of the random walk behaviour? For example, on a finite graph if some of the rates are increased, does this necessarily improve the convergence to the stationary (uniform) distribution? How does the total time spent at each vertex depend on the rates?

Consider a Cayley graph $G$ of a group $\Gamma$ generated by a finite set $S \subset \Gamma$. Write $\bp$ for the identity element and $|x|$ for the graph distance between $\bp$ and $x$. Given $r_s > 0$ for each $s \in S$, let $r(e) := r_s$ when $e$ is an edge of $G$ corresponding to $s \in S$. Throughout this paper, we will assume that generating sets are symmetric --- that is, if $s$ is a generator then so is $s^{-1}$, and that the corresponding rates $r_s$ and $r_{s^{-1}}$ are equal. 

Let $(Z_t)_{t \ge 0}$ be the corresponding random walk starting from $\bp$ and $\sigma(\rt) := \lim_{t \to\infty} \Ebig{|Z_t|}/t$ be its limiting speed. This limit exists by a well-known subadditivity argument; the subadditive ergodic theorem also shows that $\lim_{t \to\infty} |Z_t|/t = \sigma(\rt)$ a.s. Our main result for infinite graphs is the following:

\procl t.coxeter 
If $(\Gamma, S)$ is a Coxeter system, 
then $\sigma(\rt)$ is monotone increasing in $\rt$.
If $(\Gamma, S)$ is irreducible and nonelementary hyperbolic, then $\sigma(\rt)$ is strictly increasing in $\rt$.
\endprocl

See Corollary \ref{cor:coxmonotonespeed} and Theorem \ref{the:hyprates} for the proofs.
Cayley graphs that are trees are examples of Cayley graphs that arise from Coxeter systems.
The theorem does not extend to all Cayley graphs: see Example \ref{ex:mod2nfreeproduct}.

We then consider finite Cayley graphs and examine the convergence to stationarity. It is plausible that increasing some of the rates might always improve the convergence to the stationary distribution (which is the uniform distribution), but this is not true in general. 
Our main positive result on this topic is the following:

\procl t.ratemajorCox-intro
Let $(\Gamma, S)$ be a finite Coxeter system. 
Let $\rt$ and $\rt'$ be two sets of rates on $S$ with $r_s \le r'_s$ for all $s \in S$. Let $t > 0$. Denote the corresponding transition probabilities by $p_t(x, y; \rt)$ and $p_t(x, y; \rt')$. Then $p_t(o, \cdot; \rt)$ majorizes $p_t(o, \cdot; \rt')$ with inequality if $\rt \ne \rt'$.
\endprocl

This implies that increasing rates leads to distributions that are closer to the stationary distribution in many senses, including $\ell^p$ for $1 \le p \le \infty$ and relative entropy. See \rref d.major/ for the definition of majorization and \rref t.ratemajorCox/ for the proof. The edge graph of the permutohedron is an example of a Coxeter Cayley graph, where $\Gamma$ is a symmetric group.
For certain special Markov chains in discrete time, a much stronger inequality for total variation distance was proved by \cite[Theorem 8.3]{FillKahn}, which slightly extended a result of \cite{PeresWinkler}. Similar strong inequalities using what we call refresh rings hold in our context as well: see Theorems \briefref t.coxdiscretedist/ and \briefref t.coxdiscretemajor/.

Despite the impossibility of ``perfect'' couplings, Markovian couplings will be crucial to our proofs of Theorems \briefref t.coxeter/ and \briefref t.ratemajorCox-intro/.

In Section \ref{sec:stationarity}, we prove results of this kind for some other special cases, such as for every Cayley graph of an abelian group, or for any group when distance from stationarity is measured with either the $\ell^2$- or $\ell^\infty$-distance.
We then discuss examples where increasing a rate worsens the distance to stationarity, which we find 
in the Cayley graphs of groups as small as the dihedral group $D_5$ and the symmetric group $S_4$. How easy it is to find these examples depends on which $\ell^p$-distance is used to measure the distance to stationarity. For instance, we find examples in dihedral groups for $p$ very close to $4$, both above and below $4$, but are unable to find any for $p$ equal to $4$. In the symmetric group $S_5$, we do find examples for $p$ exactly equal to $4$.  
Hermon and Kozma \rref b.HerKoz/ find examples of generating sets for symmetric groups where increasing some rates by a tiny amount has the effect of increasing the mixing time by a large amount, which also means that the $\ell^1$-distance increases by a large amount. Their examples have unbounded degree as the size of the group tends to infinity.

A third result goes beyond the setting of Cayley graphs and concerns arbitrary rate functions $\rt$, but is restricted to the ray graph, i.e., the nearest-neighbor graph on $\Z^+$. Our result complements that of Karlin and McGregor given above.

\procl t.ray
With $\rt$ being an arbitrary positive rate function on the edges of the nearest-neighbor graph on $\Z^+$ and $t > 0$, $\>\Ebig{|Z_t|}$ is monotone increasing in $\rt$.
\endprocl

See Corollary \ref{cor:expdist} for the proof.
The theorem does not extend to walks on all of $\Z$: see Example \ref{ex:Zdist}.

Another aspect of the behaviour of random walk concerns the entropy: Let $h_t(x, \rt)$ be the (natural-log) entropy of $Z_t$ when started at vertex $x$. The following conjecture \cite[Conjecture 4.11]{BLS:pert} turns out to be false:

\procl g.ent
Let $G$ be a finite graph and $t > 0$. Then $|V(G)|^{-1} \sum_{x \in V(G)} h_t(x, \rt)$ is monotone increasing in $\rt$.
\endprocl

For a counterexample, let $G$ be a star with 6 vertices, all edges but one having rate 1 and the other having rate either 10 or 20. Direct calculation shows that the mean entropies at time 1 are, respectively, approximately $1.626355024$ and $1.626293845$.
\rref t.ratemajorCox-intro/ implies that the conjecture does hold for rates depending on generators in the setting of Coxeter systems, but changing the generators, even on dihedral groups, can make it fail: see Subsection \ref{sec:negative}.
However, for the purposes of \cite{BLS:pert}, it would suffice that there be a lower bound on the factor by which the mean entropy can decrease; we do not know whether this weakened statement is true.

Other monotonicity results for random walks include \rref b.RegShi/ and \rref b.McMP/, which study monotonicity in time of $p_t(x, y)/p_t(x, x)$, and \rref b.ChenSun/, which studies the expected range of symmetric random walks on $\Z^d$ when extra steps are included deterministically.

\section{Coxeter systems}

\subsection{Background on Coxeter systems}

A \defn{Coxeter system} $(W, S)$ is a group $W$ with generators $S$ defined by a presentation of the form $\pres{S}{\forall s, s' \in S\ (s s')^{m(s, s')} }$, where each $m(s,s')$ is either a positive integer or $\infty$. Each $m(s, s)$ is equal to $1$, and $m(s, s') \ge 2$ for $s \ne s'$. When $m(s, s') = \infty$, the interpretation is that no such relation is imposed. In particular, all generators are involutions. Furthermore, $m(s, s') \in \{1, 2\}$ iff $s$ and $s'$ commute. As always in this paper, we assume that $S$ is finite. It is customary to denote Coxeter groups with the letter $W$ (for ``Weyl''), and we will do that in order to make it easier to discern which of our results apply to Coxeter groups. Coxeter groups are abstractions of reflection groups; dihedral and symmetric groups are simple examples. Section 1.2 of \cite{BjBr} is devoted to these and other examples.

We always use \defn{right Cayley graphs}: the vertex set is a group $\Gamma$ and the unoriented edges are $\bigl\{\{x, xs\} \st x \in \Gamma, s \in S\bigr\}$. One could use multigraphs, but nothing would change for our questions of interest. The \defn{Cayley diagram} of a group $\Gamma$ with generators $S$ is the corresponding Cayley graph with edges labelled by the corresponding generator; if the generator is not an involution, then also the edge is given an orientation so that multiplication (on the right) by the generator maps the tail of the edge to the head.
For $\gamma \in \Gamma$, we write $L_\gamma$ for left multiplication by $\gamma$, which is an automorphism of the Cayley diagram of $\Gamma$.

If $(W,S)$ is a Coxeter system, then the Cayley graph of $W$ with respect to $S$ has many nice symmetry properties. One simple property is that the Cayley graph is bipartite, because all relations give even-length cycles. 
When $G$ is a bipartite Cayley graph and $L$ is a left multiplication that interchanges the endpoints of some edge, then we call $L$ a \defn{reflection}.
We also call the set $M$ of edges preserved by $L$ a \defn{wall}, a \defn{hyperplane}, or a \defn{mirror}.
A key property for our purposes is the following, which will allow us to use arguments akin to the well-known reflection principle for one-dimensional random walks and Brownian motion.

\begin{Lemma} \label{lem:walldef}
Let $w$ and $ws$ be two adjacent vertices in the Cayley graph of a Coxeter system, $(W,S)$. Let $M$ be the set of edges preserved by the automorphism $L := L_{wsw^{-1}}$. The map $L$ interchanges the endpoints of each edge of $M$. The wall $M$ separates the vertices closer to $w$ from those closer to $ws$. Suppose that $|w| < |ws|$; then for all $v \in W$, $v$ is closer to $w$ than to $ws$ iff $|v| < |L v|$.
\end{Lemma}

These properties of walls are well known, but not all are easy to find explicitly stated in this form. Because they are crucial to our results, we provide a proof.
Denote the shortest-path metric in a graph by $\dist$.

\rproof
Let $G$ be the Cayley graph. Refer to vertices closer to $w$ than to $ws$ as white and the others as black; also, call an edge grey if its endpoints have different colours. If $y$ is white and $z$ is black, then every path from $y$ to $z$ must include a grey edge. 

Merely because $G$ is a bipartite graph, no vertex is equidistant from $w$ and $ws$. Because $L$ preserves distances and interchanges $w$ and $ws$, it therefore also interchanges white and black vertices. In particular, every preserved edge is grey, and the endpoints of grey edges are interchanged. We need to use more than bipartiteness to show that all grey edges belong to the wall.

Let $(x, xs')$ be a grey edge with $x$ white.
First, merely by bipartiteness, $\dist(w, x) = \dist(ws, xs')$. We may write this last conclusion as $|w^{-1}x| = |sw^{-1}xs'|$. Suppose that $\wrd$ is a reduced word for $w^{-1}x$. Then $s\wrd s'$ is a longer word that, when reduced, has the same length as $\wrd$. By the deletion condition of Coxeter systems (\cite[Proposition 1.4.7]{BjBr}, \cite[Corollary 5.8]{Humphreys}, or \cite[Theorem 3.2.17]{Davis}), there are two letters from $s\wrd s'$ that can be deleted in order to obtain a reduced word for the same element. Since $x$ is white and $ws$ is black, $|sw^{-1}x| = |w^{-1}x| + 1$, and so the word $s\wrd$ must be reduced; similarly, the word $\wrd s'$ is reduced. Checking cases then reveals that the two letters to delete must be the initial $s$ and the final $s'$, which yields the word $\wrd$. In other words, $sw^{-1}xs' = w^{-1}x$. Another way to say this is that $L_{wsw^{-1}} x = xs'$. In particular, $(x, xs')$ belongs to the wall, as claimed. 

Finally, suppose that $o$ is white. If $v$ is white, then let $(x, xs') \in M$ be an edge belonging to a geodesic from $o$ to $L v$. We have $|v| \le \dist(o, x) + \dist(x, v) = \dist(o, x) + \dist(xs', L v) = |L v| - 1 < |L v|$.
If $v$ is, instead, black, then $L v$ is white, so $|L v| < |L(L v)| = |v|$.
\Qed

\procl c.reflections
Let $(W, S)$ be a Coxeter system. The set of reflections equals the set of left multiplications $L_{wsw^{-1}}$ for $w \in W$ and $s \in S$. If $a, b \in S$ and $w, x \in W$, then $(w, wa)$ and $(x, xb)$ belong to the same wall iff $waw^{-1} = xbx^{-1}$.
\endprocl

\rproof
It is clear that $L_{wsw^{-1}}$ interchanges the endpoints of the edge $(w, ws)$, and, conversely, if $L_\gamma$ interchanges those endpoints, then $\gamma w = w s$, i.e., $\gamma = wsw^{-1}$. For the last statement in the corollary, if the two edges belong to the same wall, defined, say, by a reflection $L_\gamma$, then by Lemma \ref{lem:walldef}, $\gamma w = wa$ and $\gamma x = xb$, whence $waw^{-1} = \gamma = xbx^{-1}$. The converse is proved similarly.
\Qed

Thus, there is a natural bijection between reflections and walls.
 We will identify a wall $M$ with the (disconnected) graph formed by the edges $E(M)$ in the wall and their endpoints $V(M)$.

\begin{figure}[ht]
\begin{center}
\includegraphics[width=\textwidth]{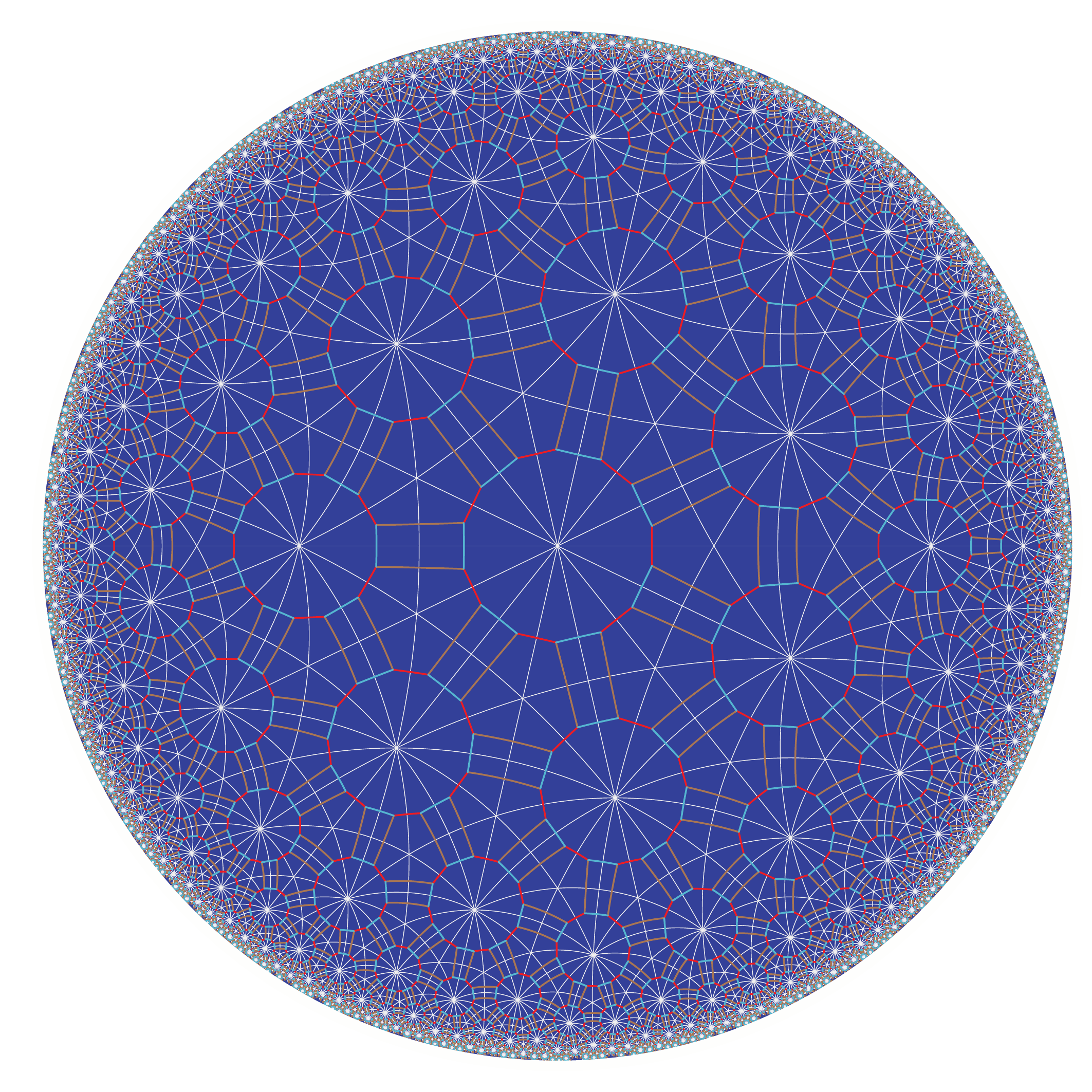}
\end{center}
\caption{A portion of the Cayley graph of the $(2, 3, 7)$-triangle group, drawn by Matthias Weber. Edges are coloured red, light blue, or light brown according to the corresponding generator. Each dark blue face corresponds to a relation. Each white geodesic corresponds to a wall.} 
		\label{fig:237cayley}
\end{figure}

In the case of a free Coxeter group (i.e., the free product of copies of $\Z/2\Z$), whose Cayley graph is a tree, walls are the same as single edges.
Figure \ref {fig:237cayley} shows the Cayley diagram of the Coxeter group $\langle a, b, c \mid a^2,b^2,c^2, (ab)^7, (bc)^2, (ca)^3 \rangle$ corresponding to reflections in the sides of a hyperbolic triangle of angles $(\pi/2, \pi/3, \pi/7)$. 
Each white geodesic crosses through the edges of a wall. Note that in this case, each wall contains edges corresponding to each of the three generators.

When we prove strict monotonicity of speed in the rates, we will need 
the following three lemmas, the first two of which are well known.  

\procl l.crossing
Let $M$ be a wall in a Coxeter Cayley graph. Then every geodesic uses at most one edge in $M$.
\endprocl

\rproof
It suffices to show that if $(w, wa), (v, vb) \in M$ both belong to a path from $w$ to $v$, then there is a shorter path from $w$ to $v$. Indeed, if $P$ is a subpath from $wa$ to $vb$, then the path obtained from $P$ by applying the reflection in $M$ is a path from $w$ to $v$ that is shorter by at least two edges. 
\Qed

\procl l.stays
Let $(W, S)$ be a Coxeter system. There is some finite $K$ such that for every wall $M$, every geodesic between two vertices of $M$ stays within distance $K$ of $M$.
\endprocl

\rproof
By the parallel-wall theorem of Brink and Howlett (\cite{BH} or \cite[Theorem C]{Caprace}), there is some $K$ so that for any pair $(w, M)$ where $w$ is a vertex at distance at least $K$ from a wall $M$, there is some other wall $M' \ne M$ that separates $w$ from $M$. By \rref l.crossing/, no geodesic between vertices of $M$ can cross $M'$, whence it cannot go farther than $K$ from $M$.
\Qed

\procl l.qtrwall
Let $M$ be a wall in a Coxeter system $(W, S)$. Let $W_M := \{w \in W \st L_w V(M) = V(M)\}$. Then $W_M$ acts quasi-transitively by graph automorphisms of $M$, in other words, there are only finitely many orbits in $V(M)$ under the left action of $W_M$. Moreover, the orbit of $w \in V(M)$ is $\{v \in V(M) \st (v, vs) \in E(M)\}$, where $s$ is the generator for which $(w, ws) \in E(M)$.
\endprocl

\rproof
Suppose that $(w, ws), (v, vs) \in E(M)$. Clearly $L_{vw^{-1}}$ carries the first edge to the second. We claim that $vw^{-1}  \in W_M$.
To prove our claim, let $(x, x a) \in E(M)$. We have $wsw^{-1} = vsv^{-1} = xax^{-1}$ by \rref c.reflections/, whence $vw^{-1}xa(vw^{-1}x)^{-1} = vsv^{-1}$, so that indeed $(vw^{-1}x, vw^{-1}xa) \in E(M)$. Therefore, for each pair of edges in $M$ corresponding to the same generator, there is an element of $W_M$ that takes one edge to the other. Since $S$ is finite, it follows that this action is quasi-transitive. Finally, since left multiplications preserve the Cayley diagram, the orbit of $w$ cannot be any larger than what is claimed.
\Qed

\subsection{Monotonicity of speed}

A key technical device we will use to construct our couplings is the following alternative way to generate our random walks.

\begin{Definition}\label{def:refresh}
If an edge $e$ has an attached Poisson clock of rate $r$, we will generally consider this process as being controlled by a Poisson clock of rate $2r$, and when this new clock rings, there is a $\frac{1}{2}$ chance that the original clock rings, otherwise nothing happens. We will reserve the upper case $R$ (for ``refresh time'') for this Poisson process of twice the rate, with subscripts denoting the edge in question, for instance, $R_e$. We may use a fair coin flip to decide whether the original clock rings; if so, we will say the coin flip is ``move'', otherwise that the coin flip is ``stay''. Alternatively, if the walk is at an endpoint of $e$ at the time that $R_e$ rings, then we may randomise the walk immediately after that ring to be at either endpoint of $e$ with equal probability. In the case of Cayley graphs where the rates depend only on the generators, we observe that it suffices to use only one Poisson process per generator rather than one per edge; we then use a subscript corresponding to the generator.
\end{Definition}

Although our goal is \rref t.coxeter/, we begin for concreteness and clarity with a specific example, the 3-regular tree.

\begin{Example}
\label{ex:treecayley}
Consider the $3$-regular tree as the Cayley graph of the group $\pres{a,b,c}{a^2,b^2,c^2}$, so that each edge is associated to one of the three generators, $a$, $b$, or $c$. Choose rates $r_a$, $r_b$ and $r_c$ for each generator, and use those as the rates for each corresponding edge. 
\end{Example}

It is well known that the random walk of Example \ref{ex:treecayley} is transient if and only if all three rates are positive.

Before we can discuss the dependence of the escape speed on the rates $r_a$, $r_b$, and $r_c$, we need some preliminary results.

\begin{Proposition}\label{prop:monotoneprobs}
Consider Example \ref{ex:treecayley}. If $w$ and $w'$ are two adjacent vertices with $w$ closer to the initial position than $w'$, then for all $t > 0$, we have $\P[Z_t = w] > \P[Z_t = w']$. 
\end{Proposition}
\begin{proof}
Without loss of generality, assume that the edge between $w$ and $w'$ is labelled by the generator $a$, so that $w' = wa$. 

We prove our result by giving a probability-preserving injection from paths (of the random walk up to time $t$) resulting in $w'$ to paths resulting in $w$. For any path that ends at $w'$, let $T$ be the first time $\le t$ at which the path is at $w$ and $R_a$ rings. Modify the path by changing the outcome of the coin flip attached to this ring of $R_a$, and proceed with the rest of the path by multiplying by $a$, $b$, or $c$ when appropriate, not changing any of the other rings of the Poisson processes or the coin flips. This has the effect of applying the reflection $L_{waw^{-1}}$ to each state after time $T$, so paths that ended at $w'$ now end at $w$.

Indeed, this map is a probability-preserving bijection from paths that ever cross the edge from $w$ to $w'$ and end at $w'$ to paths that end at $w$ and are at $w$ at some ring of $R_a$. Crucially, it is impossible to get to $w'$ without crossing the edge $(w, w')$, and hence without being at $w$ when $R_a$ rings. However, there is a positive probability that the walk ends at $w$ without ever being at $w$ when $R_a$ rings, giving our strict inequality.
\end{proof}

Proposition \ref{prop:monotoneprobs} actually has a much shorter proof: 
Let $\tau := \inf\{t > 0 \st Z_t = w\}$, and let $p_t(x, y)$ denote the transition probabilities of $Z_t$. Then
\eqaln{
\P[Z_t = w']
&=
\Ebig{\P[Z_t = w' \mid \tau]}
=
\E[p_{t - \tau}(w, w') \st \tau < t]
\\ &<
\E[p_{t - \tau}(w, w) \st \tau < t]
=
\P[Z_t = w]
}
(compare \eqref{eq:infty}).
The reason we did not give this proof is that
Proposition \ref{prop:monotoneprobs} may be generalised to the Cayley graph of any Coxeter system, but not with this short proof. This generalisation is a continuous-time version of Theorem 1 of \cite{Bruhat}, and the proof is essentially the same.
If $M$ is a wall, let us call a time $t$ an \defn{$M$-refresh time} for our random walk if $Z_t \in M$ and the refresh clock $R_a$ rings at time $t$, where $a$ is the generator such that $(Z_t, Z_t a) \in M$. 
The main difference between the case of trees and the general case of Coxeter systems is that in the tree, every path that reaches $w'$ must cross the edge between $w$ and $w'$. In a general Coxeter system, this is no longer the case, so instead we need to consider the steps at which the walk might move from a state closer to $w$ to one closer to $w'$, in other words, the $M$-refresh times for the wall determined by $(w, w')$.

\begin{Theorem}\label{the:coxmonotoneprobs}
Let $(W, S)$ be a Coxeter system. To each generator $s \in S$, attach a Poisson clock of rate $r_s$. Consider a random walk starting at the identity that moves from a current location $w$ to $ws$ when the $s$-clock rings. If $w$ and $w'$ are two adjacent vertices with $|w| < |w'|$, then for all $t > 0$, we have $\P[Z_t = w] > \P[Z_t = w']$. 
\end{Theorem}
\begin{proof}
Let $w' = wa$. We follow the proof of Proposition \ref{prop:monotoneprobs} in this more general setting.

Fix $t > 0$.
We construct a probability-preserving injection from paths of the random walk resulting in $w'$ to paths resulting in $w$. Let $M$ be the wall determined by $(w, w')$ and $L_M$ be the reflection in $M$. Let $T$ be the first $M$-refresh time. 

For a path ending at $w'$ at time $t$, we have $T < t$.
Modify such a path by changing the outcome of the coin flip attached to the refresh ring at time $T$, and proceed with the rest of the path by leaving the remaining sequence of refresh rings and coin flips unchanged. This has the effect of applying the reflection $L_M$ to each state after time $T$, so paths that ended at $w'$ now end at $w$. 

This map is a probability-preserving bijection from the set of paths that end at $w'$ at time $t$ to the set of paths that end at $w$ at time $t$ and that have $T < t$.
There is a positive probability that the walk ends at $w$ without taking such a path,
giving our strict inequality.
\end{proof}

\begin{Remark}\label{rem:condition1}
Let $(u_i)$ be any sequence of times and $(a_i)$ be any sequence of generators. Note that Theorem \ref{the:coxmonotoneprobs} remains true if the random walk is conditioned to have $R_{a_i}$ ring at time $u_i$ for each $i$. The same holds if, in addition, we condition that $R_s$ has no other rings for certain generators $s$ up to time $t$. (Critically, $R_{a_i}$ is the refresh clock of Definition \ref{def:refresh}, not the clock that always multiplies by the generator $a_i$.)

The proof is unchanged, because the constructed injection preserves the times at which each $R_s$ rings.
\end{Remark}

We can express our argument in another way. The location $Z_t$ is a function of the times and generators of the refresh rings that occur before time $t$ and the results of the corresponding coin flips. Let a coin flip 0 represent ``stay'' and 1 represent ``move''.
For $s \in S$, let $\cfm(s,0)$ be the identity element and $\cfm(s,1) := s$. Given sequences $\seqs = (s_1, \ldots, s_n) \in S^n$ and $\bits = (b_1, \ldots, b_n) \in \{0, 1\}^n$, let 
\[
\cfm(\seqs, \bits) := \cfm(s_1,b_1)\cfm(s_2,b_2) \cdots \cfm(s_n,b_n) \in W.
\]
Thus, if $\seqs$ gives the sequence of generators whose refresh rings occur before time $t$ and $\bits$ is the corresponding sequence of coin flips, we have $Z_t = \cfm(\seqs, \bits)$.
For any sequence $(a_1, \ldots, a_n)$ and $0 \le k \le n$, write $(a_1, \ldots, a_n)_k$ for the initial segment $(a_1, \ldots, a_k)$. Our argument shows the following, which also easily implies Theorem \ref{the:coxmonotoneprobs}:

\procl t.coxdiscreteprobs
Let $(W, S)$ be a Coxeter system. Let $\seqs = (s_1, \ldots, s_n)$ be a finite sequence from $S$. 
Let $\cfs = (\cf_1, \ldots, \cf_n)$ be independent, uniform $\{0, 1\}$-valued random variables.
If $w$ and $w'$ are two adjacent vertices with $|w| < |w'|$, then $\P[\cfm(\seqs, \cfs) = w] \ge \P[\cfm(\seqs, \cfs) = w']$. 
\endprocl

\rproof
Let $M$ be the wall determined by $(w, w')$ and $L_M$ be the reflection in $M$. Given $\bits = (b_1, \ldots, b_n)$, define $\bits' := (b_1, \ldots, b_j, 1 - b_{j+1}, b_{j+2}, \ldots, b_n)$ if there is some $k < n$ such that $\cfm(\seqs_k, \bits_k) \in M$ and $j$ is the smallest such index, while $\bits' := \bits$ if there is no such $k$. Note that $(\bits')' = \bits$. Thus, $\bits \mapsto \bits'$ is a permutation of $\{0, 1\}^n$ with the property that $\cfm(\seqs, \bits') = L_M \cfm(\seqs, \bits)$ iff $\bits' \ne \bits$. In particular, $\cfm(\seqs, \cfs) = w'$ implies $\cfm(\seqs, \cfs') = w$. Since $\cfs'$ has the same distribution as $\cfs$, the inequality follows.
\Qed

We will now construct a Cayley graph where the conclusion of Theorem \ref{the:coxmonotoneprobs} fails. 

\begin{Example}\label{ex:mod8}
Consider the cyclic group $\Z / 8\Z$, with generators $\pm1$, $\pm2$, and $\pm3$. Let the rates be chosen so that $r_1 = r_3$ are small and $r_2 = 1$. At time $2$, there is a reasonable chance that the $2$-clock has rung exactly twice and the $-2$-clock has not rung, and the $\pm1$-clocks and $\pm3$-clocks are unlikely to have yet rung. Then the walker, started at $0$, is more likely to be at $4$ than at $3$, even though $4$ is distance $2$ from $0$, $\>3$ is distance $1$, and there is an edge between $3$ and $4$. Thus, Proposition \ref{prop:monotoneprobs} is not true in this setting.

The issue here is that while the vertices $3$ and $4$ are connected by an edge, and $4$ is farther from $0$ than $3$ is, we have made this edge less likely
than others, and there are paths from $0$ to $4$ that do not need such an unlikely edge. In the Coxeter-system setting, this cannot happen according to Theorem \ref{the:coxmonotoneprobs}.
\end{Example}

We may also modify this example so that the expected distance is not monotone in time; a discrete-time example of this phenomenon is due to Oded Schramm \cite[Exercise 13.16(b)]{LP:book}.

\begin{Example}\label{ex:mod2n}
Adjust Example \ref{ex:mod8} so that the group is $\Z / 2n\Z$, with $n$ at least $5$. Take as generators $\pm 2$ and every odd number, with $\pm 2$ having high rates and the odd numbers having low rates.

After a time chosen so that the $\pm 2$-clocks have rung many times and the odd generator-clocks are unlikely to have rung, the walker is close to uniformly distributed on the even numbers, with expected distance from $0$ close to $2 - \frac{4}{n}$. After much more time, the walker will be close to uniformly distributed on all $2n$ states, with expected distance close to $\frac32 - \frac{2}{n}$. Thus, the expected distance is not monotone.
\end{Example}

This allows us to produce an example where the escape speed is not monotone in the rates.

\begin{Example}\label{ex:mod2nfreeproduct}
Take the free product of Example \ref{ex:mod2n} with $\Z/2\Z$, using the same generators as in that example, and one new generator $a$ for
the new factor of $\Z/2\Z$. Then the escape speed is not monotone in the rates $r_2 = r_{-2}$ and $\ro$ associated, respectively, to $\pm 2$ and the odd generators of the $\Z / 2n\Z$ factor.
\end{Example}
\begin{proof}
Take the rate $r_a$ to be $1$, $r_2$ to be very large, $\ro$ to be $0$, and $n$ to be large. Then the walker moves through copies of Example \ref{ex:mod2n} at rate $1$, backtracking only with rate proportional to $\frac{1}{n}$, and the average distance between entry and exit points in each copy of $\Z / 2n\Z$ is $2 - \frac{4}{n} - o_{r_2}(1)$. Therefore the escape speed is $3 - O(\frac{1}{n}) - o_{r_2}(1)$. 

Now, increase the rate $\ro$ to be as large as $r_2$. The average distance between entry and exit points of each copy of $\Z / 2n\Z$ is now only $\frac32 - \frac{2}{n} - o_{r_2}(1)$, so the escape speed has decreased to $\frac52 - O(\frac{1}{n}) - o_{r_2}(1)$.
\end{proof}

\procl r.separate-monot
This example can be modified so that the escape speed is not monotone, either increasing or decreasing, in each rate separately, even for a rate on an involution.
%
\endprocl

The result of Example \ref{ex:mod2nfreeproduct} should not necessarily be surprising --- the graph distance is calculated using all edges, regardless of associated rates, so we shouldn't expect this necessarily to be a terribly meaningful quantity when it is affected by edges of very low rates. The following results start with trees, where this cannot occur, and we then show that Cayley graphs of Coxeter systems are symmetric enough that our results still apply. (While Examples \ref{ex:mod2n} and \ref{ex:mod2nfreeproduct} could equally well be implemented in Coxeter groups, replacing finite cyclic groups and $\Z$ with finite and infinite dihedral groups, the generators required would not be the Coxeter generators. Our results for Coxeter groups require that the generating set is composed of the Coxeter generators.)

We now return to settings Example \ref{ex:treecayley}.

\begin{Proposition}\label{prop:extrars}
Let $t_1$ be an arbitrary time, and let $s \in S$ be one of the generators. Consider the following two random processes:
\begin{enumerate}
\item Run an instance $Z^1$ of the random walk of Example \ref{ex:treecayley}.
\item Run a separate instance $Z^2$ of the random walk of Example \ref{ex:treecayley} until time $t_1$, then trigger $R_s$ at the (deterministic) time $t_1$, and then continue to run the random walk.
\end{enumerate} 
These two processes may be coupled so that at all times $t > t_1$, either they are in the same state or the second walker is farther from the origin than the first is. 
\end{Proposition}
\begin{proof}
Given $Z^1$, we will construct $Z^2$ via a special coupling to $Z^1$.
Let $A := \{ w \st |w| < |ws|\}$. 
Define $\gamma$ to be the random vertex such that $Z^1_{t_1} \in \{\gamma, \gamma s\}$ with $\gamma \in A$. Couple $Z^1$ and $Z^2$ up to time $t_1$ to be independent conditional that $Z^2_{t_1} \in \{\gamma, \gamma s\}$, which is possible because $Z^1$ and $Z^2$ have the same law up to time $t_1$. 
For all $w \in A$, Proposition \ref{prop:monotoneprobs} gives us that 
\rlabel e.gamma1
{\P[Z^1_{t_1} = \gamma \mid \gamma = w]
=
\frac{\P[Z^1_{t_1} = w]} {\P[Z^1_{t_1} = w] + \P[Z^1_{t_1} = ws]} 
>
\frac12.}
We need $Z^2_{t_1^+}$ to be uniformly distributed in $\{\gamma, \gamma s\}$ independently of $(Z^2_{t})_{0 \le t \le t_1}$ conditional on $\gamma$.
We may ensure this by letting
$Z^2_{t_1^+}$ equal $\gamma s$ whenever $Z^1_{t_1} = \gamma s$, and sometimes even when $Z^1_{t_1} = \gamma$; this is possible in light of \rref e.gamma1/ and maintains the required independence by the first step in our construction of $Z^2$.
We have now achieved that $|Z^1_{t_1}| \le |Z^2_{t_1^+}|$.

In case $|Z^1_{t_1}| = |Z^2_{t_1^+}|$ (which has probability $2\P[Z^1_{t_1} = \gamma s]$), we couple the two processes so that they stay together at all times after $t_1$. Otherwise, we couple them after time $t_1$ so that the second is the first reflected by $L_{\gamma s\gamma ^{-1}}$ (which is true at time $t_1^+$), until and unless they would move across the edge between $\gamma $ and $\gamma s$, in which case we couple the $R_s$-coin flips so that the walks agree from then on. 

For any time $t > t_1$, this results in either the two processes being at the same vertex at time $t$, or the first being at a vertex $v$ that is closer to $\gamma $ than to $\gamma s$, and the second being at $L_{\gamma s\gamma ^{-1}}(v)$. In this case, $v$ is closer to the initial state than $L_{\gamma s\gamma ^{-1}}(v)$ is by Lemma \ref{lem:walldef}. 
\end{proof}

\begin{Remark}\label{rem:condition2}
As with Remark \ref{rem:condition1}, let $(u_i)$ be any sequence of times in $[0, \infty) \setminus \{t_1\}$ and $(a_i)$ be any sequence of generators. Proposition \ref{prop:extrars} remains true if the processes $Z^1$ and $Z^2$ are each conditioned to have $R_{a_i}$ ring at time $u_i$ for each $i$.
The same holds if, in addition, we condition that $R_s$ has no other rings up to any given time $t_2$ for the process $Z^1$. 

The proof is unchanged, except for applying Remark \ref{rem:condition1} when we invoke Proposition \ref{prop:monotoneprobs}.
\end{Remark}

\begin{Corollary}\label{cor:increasingspeed}
The escape speed in Example \ref{ex:treecayley} is nondecreasing in the rates $r_a$, $r_b$, and $r_c$.
\end{Corollary}
\begin{proof}
By symmetry, it suffices to show that the escape speed is nondecreasing in $r_a$. It suffices for this to show that $\Ebig{|Z_t|}$ is nondecreasing in $r_a$ for every time $t$. But increasing the rate $r_a$ just results in extra instances of $R_a$, so this is a consequence of Proposition \ref{prop:extrars} as follows.

Let $(X_t)$ be our random walk and $(Y_t)$ be the random walk where the rate $r_a$ has been increased. We may couple $X$ and $Y$ so that they have the same $R_a$, $R_b$, and $R_c$ rings at the same times, except that $Y$ has some random number $N_t$ of additional rings of $R_a$ up to time $t$. For each $i$ between $0$ and $N_t$, let $X^i$ be the random walk which agrees with $X$ except that it has the first $i$ of these additional rings of $R_a$, so that $X^0 = X$ and $X^N = Y$. 

For each $i$, we may use Proposition \ref{prop:extrars} to couple $X^i$ and $X^{i+1}$ so that at time $t$ the walk $X^{i+1}$ is at least as far from the origin as $X^i$ is. Because $X^{i+1}$ has $R_a$-rings at the times that $X^i$ does, this application of Proposition \ref{prop:extrars} requires the observation of Remark \ref{rem:condition2}; we are also conditioning on the time of the extra $R_a$-ring. Combining these couplings gives a coupling between $X = X^0$ and $X^N = Y$ where at time $t$ the walk $Y$ is at least as far from the origin as $X$ is.
\end{proof}

As with Proposition \ref{prop:monotoneprobs}, these results also apply to walks on the Cayley graphs of Coxeter systems.

\begin{Corollary}\label{cor:coxmonotonespeed} 
The result of Proposition \ref{prop:extrars} is true for the Cayley graph of a Coxeter system, as is Corollary \ref{cor:increasingspeed}.
\end{Corollary}
\begin{proof}
This result follows from Theorem \ref{the:coxmonotoneprobs} in the same way as Proposition \ref{prop:extrars} and Corollary \ref{cor:increasingspeed} follow from Proposition \ref{prop:monotoneprobs}, with a single change:
in the proof of Proposition \ref{prop:extrars}, the phrase

\begin{center}``until and unless they would move across the edge between $\gamma$ and $\gamma s$, in which case we couple the $R_s$-coin flips \dots ''\end{center}
should be replaced by 

\begin{center}``until and unless they would move across an edge $(x, xb)$ in the wall determined by $(\gamma, \gamma s)$, in which case we couple the $R_b$-coin flips \dots ''. \qedhere\end{center}
\end{proof}

\procl r.semicox
A superficially more general result is that Corollary \ref{cor:coxmonotonespeed} holds for Cayley graphs of presentations 
\[
\Gamma := \pres{S_1, S_2}{\forall s_1 \in S_1\  \forall s \in S_1 \cup S_2\enspace (s_1 s)^{n(s_1, s)},\ \forall s_2 \in S_2\enspace s_2^{2n(s_2)} },
\]
where each $n(s_1,s)$ and $n(s_2)$ is either a positive integer or $\infty$, each $n(s_1, s_1)$ is equal to $1$, and $n(s_1, s) \ge 2$ for $s_1 \ne s$. To see why Corollary \ref{cor:coxmonotonespeed} holds for such $\Gamma$, let $\bar S_2$ be a copy of $S_2$ and write $s_2 \mapsto \bar s_2$ for a bijection from $S_2 \to \bar S_2$. Consider the Coxeter system $W := \pres{S_1, S_2, \bar S_2}{\forall s, s' \in S_1 \cup S_2 \cup \bar S_2 \enspace (s s')^{m(s, s')}}$ with $m(s_1, s_2) := n(s_1, s_2) =: m(s_1, \bar s_2)$ when $s_1 \in S_1$ and $s_2 \in S_2$, $m(s_1, s'_1) := n(s_1, s'_1)$ when $s_1, s'_1 \in S_1$, and, finally, $m(s_2, s_2) = m(\bar s_2, \bar s_2) = 1$ and $m(s_2, \bar s_2) := n(s_2)$ when $s_2 \in S_2$. The Cayley graphs of $\Gamma$ and $W$ are the same, but to get the unoriented Cayley diagram of $\Gamma$ from that of $W$, replace each label $\bar s_2$ with the label $s_2$. As long as the generator rates for a random walk on $W$ have the property that the rate for $\bar s_2$ is the same as the rate for $s_2$, then the random walks on the two diagrams have the same law. The simplest case is when $S_1 = \varnothing$ and $n \equiv \infty$ on $S_2$, in which case $\Gamma$ is a free group with free generators.
\endprocl

Similarly to \rref t.coxdiscreteprobs/, we may deduce from our arguments the following, which also easily implies Corollary \ref{cor:coxmonotonespeed}:

\procl t.coxdiscretedist
Let $(W, S)$ be a Coxeter system. Let $\seqs$ and $\seqs'$ be finite sequences from $S$ of lengths $n$ and $n'$, respectively, with $\seqs$ a proper subsequence of $\seqs'$. 
Let $\cfs$ be a Bernoulli$(1/2)$ process.
Then $\Ebig{|\cfm(\seqs', \cfs_{n'})|} \ge \Ebig{|\cfm(\seqs, \cfs_{n})|}$. 
\qed
\endprocl

\subsection{Strict monotonicity of speed}

Corollaries \ref{cor:increasingspeed} and \ref{cor:coxmonotonespeed} show that the escape speed is an increasing function of the rates $r_s$, for each generator $s$. We now consider the question of whether this function is \emph{strictly} increasing. First, we give an example where this is not the case.

\begin{Example}\label{ex:ref:thicktree}
Consider the Coxeter system $$\pres{a,b,c}{a^2,b^2,c^2} \times \pres{d}{d^2}.$$ This group is the direct product of a free Coxeter group, whose Cayley graph is the $3$-regular tree, and the two-element group. Increasing the rate $r_d$ does not change the escape speed of the random walk. As long as the three rates $r_a$, $r_b$, and $r_c$ are all positive, increasing any of them does increase the escape speed, but the escape speed is zero if one of those three rates is zero.
\end{Example}

More generally, we have the following behaviour in product groups:

\begin{Proposition}\label{prop:productrates}
Suppose that $\Gamma_i$ is generated by $S_i$ and has identity element $\bp_i$ for $i = 1, 2$, and let $\rt_i \colon S_i \to [0, \infty)$. Then for $\Gamma := \Gamma_1 \times \Gamma_2$ generated by $\bigl(S_1 \times \{\bp_2\}\bigr) \cup \bigl(\{\bp_1\} \times S_2 \bigr)$ with rates $\rt_1 \cup \rt_2$,
the escape speed in $\Gamma$ is the sum of the escape speeds in $\Gamma_i$.
\end{Proposition}
\begin{proof}
We have $|(v_1, v_2)| = |v_1| + |v_2|$ for $v_i \in \Gamma_i$. In addition, if we write the random walk as $Z_t = (Z_t^{(1)}, Z_t^{(2)})$, then $Z_t^{(i)}$ are independent random walks in $\Gamma_i$ with rates $\rt_i$.
\end{proof}

The \defn{Coxeter diagram} of a Coxeter system $(W, S)$ is the unoriented graph with vertices $S$, edges $(s, s')$ for $m(s, s') \ge 3$, and edge labels $m(s, s')$ when $m(s, s') \ge 4$.
A Coxeter system is called \defn{irreducible} if its Coxeter diagram is connected. Every Coxeter system can be expressed as a direct product of irreducible Coxeter systems, namely, the subgroups generated by the vertices in the connected components of its Coxeter diagram. An infinite, irreducible Coxeter group is either \defn{affine}, being a finite extension of an abelian group, or \defn{nonelementary hyperbolic}, being nonelementary word hyperbolic: see \cite[Chapter 12, especially Section 12.6]{Davis} and \cite[Sections 6.8 and 6.9]{Humphreys} for more detailed information on hyperbolic Coxeter groups.

We quickly review the basic definitions for hyperbolic metric spaces. A metric space $X$ in which every pair of points can be joined by a geodesic path is called \defn{$\delta$-hyperbolic} if for each geodesic triangle in $X$, each side lies within the closed $\delta$-neighbourhood of the union of the other two sides. If $X$ is $\delta$-hyperbolic for some $\delta \ge 0$, then $X$ is called \defn{hyperbolic}. A group is \defn{word hyperbolic} if some (hence every) Cayley graph is hyperbolic in the shortest-path metric. A word hyperbolic group is called \defn{elementary} if it is finite or virtually cyclic.

\begin{Remark}\label{rem:finaffrates}
In a finite or affine Coxeter group, the escape speed is zero regardless of the rates of each generator, so increasing the rate of a generator does not increase the escape speed. Indeed, the speed is zero of every nearest-neighbor, reversible, discrete-time random walk on a graph of subexponential volume growth and bounded edge weights, a consequence of the Varopoulos--Carne inequality \cite[Theorem 13.4]{LP:book}.
\end{Remark}

Thus, the same behaviour as in Example \ref{ex:ref:thicktree} appears for any direct product of a nonelementary hyperbolic Coxeter group with a finite or affine Coxeter group. 

We will now prove that in an irreducible, nonelementary hyperbolic Coxeter system, increasing the rate of any generator strictly increases the escape speed. With the preceding remark, this will characterise those generators where increasing the rate must strictly increase the escape speed.

\begin{Theorem}\label{the:strictgenerators}
Let $W$ be a Coxeter system. Write $W$ as a direct product of irreducible Coxeter systems. If the rate corresponding to each generator is positive, then increasing a rate strictly increases the escape speed if and only if that generator belongs to a nonelementary hyperbolic irreducible factor, rather than a finite or affine factor. 
\end{Theorem}

To prove this, it suffices to examine random walks on irreducible, nonelementary hyperbolic Coxeter systems, because of Proposition \ref{prop:productrates} and Remark \ref{rem:finaffrates}.

We will need the following facts about the boundary of a Coxeter group.
We begin by reviewing the definition of Gromov boundary and some of its basic properties.
Recall that a graph is equipped with its shortest-path metric, $\dist$.
A \defn{geodesic ray} in a graph $G$ is a semi-infinite path all of whose finite subpaths are geodesic. Two geodesic rays $(x_n)_n$ and $(y_n)_n$ are \defn{asymptotic} if $\sup_n \dist(x_n, y_n) < \infty$. This is an equivalence relation, whose equivalence classes form the \defn{boundary} $\bd G$ of $G$. If $\xi \in \bd G$ and $x \in V(G)$, a geodesic ray starting from $x$ that belongs to $\xi$ will also be referred to as a geodesic from $x$ to $\xi$. We define a topology on $V(G) \cup \bd G$ as follows. Fix a vertex $o$ of $G$. Given $x_k, y \in V(G) \cup \bd G$ for $k \in \N$, we say that $\lim_{k \to\infty} x_k = y$ if for some $n_\infty, n_k \in \N \cup \{\infty\}$ ($k \in \N$), there are geodesics $(x_{k, n})_{n < n_k}$ from $o$ to $x_k$ and a geodesic $(y_n)_{n < n_\infty}$ from $o$ to $y$ such that for every $n < n_\infty$, $\> x_{k, n} = y_n$ for all sufficiently large $k$. The closed sets are then those that are closed under sequential convergence. 
In case $G$ is hyperbolic, $\bd G$ is usually referred to as its \defn{Gromov boundary}.

Consider now a Cayley graph with associated continuous-time Markov chain $(Z_t)_t$.
The embedded discrete-time random walk records only the changes made by $(Z_t)_t$. This random walk has transition probability from $w$ to $ws$ equal to $r_s/\sum_{a \in S} r_a$, where $S$ is the generating set.
We can recover the law of $(Z_t)_t$ from its embedded discrete-time random walk by jumping at the times of a Poisson process of rate $\sum_{a \in S} r_a$.
Many results for $(Z_t)_t$ follow from their corresponding results for the embedded random walk, such as this:

\begin{Lemma}[Theorem 7.4 of \cite{Kaimanovich}]
\label{proposition:convergence}
Random walk on a nonelementary word-hyperbolic group converges a.s.\ to a
boundary point, provided all rates are positive. \qed
\end{Lemma}

The law on $\bd G$ of $\lim_{t \to\infty} Z_t$ is called \defn{harmonic measure}, which depends on the starting point, $Z_0$.

The following result is due to Gromov \cite{Gromov}:
\begin{Lemma}[Lemma 2.4.4 of \cite{Calegari}] 
\label{proposition:minimal}
The action of a nonelementary word-hyperbolic group on its Gromov boundary is minimal, i.e., every orbit is dense. \qed
\end{Lemma}

\procl l.support
The harmonic measure of random walk on a nonelementary word-hyperbolic group has full support in its Gromov boundary. 
\endprocl

\rproof
Suppose that $\xi$ is a point of the support when the walk starts at $o$. Then $\gamma \xi$ is a point of the support when the walk starts at $\gamma$. Since the random walk starting at $o$ has positive chance to reach $\gamma$, it follows that $\gamma \xi$ is also in the support when the random walk starts at $o$. Thus, the result follows from Lemmas \ref{proposition:convergence} and \ref{proposition:minimal}. 
\Qed

We will use the following result of \cite{SpeyerOverflow}.

\procl l.icc
The conjugacy class of a reflection in an infinite, irreducible Coxeter system $(W, S)$ is always infinite.
\endprocl

\rproof 
Let $S$ be ordered as $(s_1, s_2, \ldots, s_r)$, where $r := |S|$. Speyer \cite{Speyer} showed that $s_1 s_2 \cdots s_r s_1 s_2 \cdots s_r \cdots s_1 s_2 \cdots s_r$ is a reduced word for any number of repetitions of $s_1 s_2 \cdots s_r$. Write $w := s_1 s_2 \cdots s_r \in W$. Combining Speyer's result with \rref l.crossing/, we deduce that the walls corresponding to the edges $(w^k, w^k s_1)$ are all distinct, whence the elements $w^k s_1 w^{-k}$ are all distinct. Since the ordering of $S$ was arbitrary, every generator has infinitely many conjugates. Therefore, so does every conjugate of a generator, i.e., every reflection.
\Qed

\begin{Lemma} \label{proposition:offwall}
Let $M$ be a wall in an irreducible Coxeter system with Cayley graph $G$. Then the closure $\overline M$ does not include all of $\bd G$.
\end{Lemma}

\rproof
Let $(v, vs) \in M$.
By \rref l.icc/, there exist $w_k \in W$ such that $|w_k^{-1} vsv^{-1} w_k| \allowbreak \to\infty$ as $k \to\infty$.
By Lemma \ref{lem:walldef}, for every $w \in W$, 
\[
| w^{-1}vsv^{-1}w| = \dist(w, L_{v s v^{-1}} w) = 1 + 2\, \dist\bigl(w, V(M)\bigr).
\]
Therefore, $\dist\bigl(w_k, V(M)\bigr) \to\infty$ as $k \to\infty$.
Let $v_k \in V(M)$ satisfy $\dist(w_k, v_k) = \dist\bigl(w_k, V(M)\bigr)$. By taking a subsequence if necessary, we may assume that there is some fixed $a \in S$ such that $(v_k, v_k a) \in E(M)$ for all $k$. If we translate the geodesic that goes from $v_k$ to $w_k$ so that it starts at $v_1$, then as $k \to\infty$, it has a limit geodesic ray $(y_n)_n$, up to taking another subsequence if necessary. By \rref l.qtrwall/, $M$ is fixed under this translation, so $(y_n)_n$ is a geodesic ray from $v_1$ with $\dist\bigl(y_n, V(M)\bigr) = \dist(y_n, v_1) = n$. By \rref l.stays/, it follows that the equivalence class of $(y_n)_n$ does not belong to $\overline M$, as desired.
\Qed

\procl l.walltransient
Let $M$ be a wall in an irreducible, nonelementary hyperbolic Coxeter system with Cayley graph $G$. Provided that all rates are positive, there is some $v \in V(M)$ such that the probability that random walk started at $v$ never returns to $M$ is positive.
\endprocl

\rproof
If not, then by \rref l.qtrwall/ and the strong Markov property, the number of visits to $M$ would be infinite a.s.\ given that it starts on $M$. In other words, if the random walk starts on $M$, then it has a limit in $\overline M \cap \bd G$ a.s.\ by Lemma \ref{proposition:convergence}. But this contradicts Lemmas \ref {Lemma:support} and \ref{proposition:offwall}. 
\Qed

Given a Coxeter system $(W, S)$ and rates $\rt$, let $\eta(W, S; \rt)$ be the minimum over all walls $M$ of the probability that there is no positive $M$-refresh time for the random walk with rates $\rt$.

\procl l.norefresh
Let $(W, S)$ be an irreducible, nonelementary hyperbolic Coxeter system and $\rt$ have all positive rates. Then $\eta(W, S; \rt) > 0$.
\endprocl

\rproof
Given any wall, $M$, there is an automorphism of the Cayley diagram that takes $M$ to a wall containing $o$. By Lemmas \briefref l.walltransient/ and \briefref l.qtrwall/, it follows that there are some $d < \infty$ and $c > 0$ such that for all walls $M$ and all vertices $w \in M$, there is a path of length at most $d$ from $w$ to some $v \notin M$ from which the probability that the random walk with rates $\rt$ never visits $M$ is at least $c$.

Given any wall, $M$, let $w$ be the first vertex in $M$, if any, visited by our random walk. 
Consider a fixed path $(w = x_0, x_1, \ldots, x_\ell = v)$ of length $\ell \le d$ from $w$ to some $v \notin M$ from which the probability that the random walk with rates $\rt$ never visits $M$ is at least $c$. By reflecting $v$ in $M$ if needed, we may assume that $w$ and $v$ are on the same side of $M$, and hence that our fixed path does not cross $M$. The chance that the next $\ell$ steps of our random walk are exactly along this path and that there is no $M$-refresh time before visiting $v$ is at least the chance that $R_{x_{i-1}^{-1} x_i}$ ($1 \le i \le \ell$) are the next $\ell$ refresh rings and that the first $\ell$ coin flips all tell the random walk to move, which equals
\[
2^{-\ell} \prod_{i=1}^\ell \frac{r_{x_{i-1}^{-1} x_i}}{\sum_{s \in S} r_s}.
\]
It follows that 
\[
\eta(W, S; \rt)
\ge
\biggl(\frac{\min_{s \in S} r_s}{2\sum_{s \in S} r_s}\biggr)^d c
> 0.
\qedhere
\]
\Qed

Gou\"{e}zel \rref b.Gou/ showed that $\rt \mapsto \sigma(\rt)$ is analytic for arbitrary generators of nonelementary word-hyperbolic groups, provided all rates in $\rt$ are strictly positive.  

\begin{Theorem}\label{the:hyprates}
If $(W, S)$ is an irreducible, nonelementary hyperbolic Coxeter system, then increasing the rate $r_s$ of any generator $s$ strictly increases the escape speed, provided all rates are positive. Moreover, $\partial \sigma(\rt)/\partial r_s \ge 2\eta(W, S; \rt)^2$.
\end{Theorem}

\rproof
In order to understand the consequences of increasing the rate of any generator, consider the construction in the proof of Corollary \ref{cor:coxmonotonespeed}, taken from the proof of Corollary \ref{cor:increasingspeed}. In that construction, we produce a sequence $X^i$ of random walk paths, each successive path having one additional ring of the refresh clock, $R_s$. The key property of that construction is that for each $i$, the paths $X^i$ and $X^{i+1}$ either eventually agree, or have $|X^{i+1}| > |X^i|$ for all large times. To show that increasing a rate strictly increases the escape speed, it suffices to show that asymptotically, a positive fraction of $i$ fall into the latter case. In our proof of this, we will be very careful to avoid dealing with dependencies among these events for different $i$, as well as with the dependencies on the number of such additional rings up to time $t$.

To be more precise, consider two sets of positive rates, $\rt$ and $\rt'$, where $\rt'$ agrees with $\rt$ except that $r'_s = r_s + \epsilon$, where $\epsilon > 0$. Let $(N_t)_t$ be a Poisson process of rate $2\epsilon$, which we use for the extra rings of $R_s$. The random walk corresponding to $\rt$ is $X = X^0$, and the random walk corresponding to $\rt$ plus $i$ additional rings of $R_s$ is $X^i$. 
We couple all $X^i$ as in the proof of Corollary \ref{cor:increasingspeed}. 
We have 
\[
\sigma(\rt') - \sigma(\rt)
=
\lim_{t \to\infty} \frac1t \sum_{i=0}^{N_t-1} \bigl(|X^{i+1}_t| - |X^i_t|\bigr) \quad\textrm{a.s.}
\]
Let $A_i$ be the event that $|X^{i+1}| > |X^i|$ for all large times. 
If $i < N_t$, then we have $|X^{i+1}_t| \ge |X^i_t|$, with equality only if $A_i$ does not occur. Thus, $|X^{i+1}_t| - |X^i_t| \ge \I{A_i}$,
whence
\eqaln{
\sigma(\rt') - \sigma(\rt)
&\ge
\limsup_{t \to\infty} \frac1t \sum_{i=0}^{N_t-1} \I{A_i}
=
\lim_{t \to\infty} \frac{N_t}t \cdot \limsup_{t \to\infty} \frac1{N_t} \sum_{i=0}^{N_t-1} \I{A_i}
\\ &=
2\epsilon \limsup_{n \to\infty} \frac1{n} \sum_{i=0}^{n-1} \I{A_i} \quad\textrm{a.s.}
}
Taking expectation and using Fatou's lemma yields
\[
\sigma(\rt') - \sigma(\rt)
\ge
2\epsilon \limsup_{n \to\infty} \frac1{n} \sum_{i=0}^{n-1} \P(A_i).
\]

Let $T_i$ be the time of the $i$th ring of $(N_t)_t$, and let $M_i$ be the wall containing the edge $(X^i_{T_i}, X^i_{T_i} s)$. Recall that a time $t$ is an $M_i$-refresh time for $X^i$ if $X^i_t \in M_i$ and the refresh clock $R_a$ rings at time $t$, where $a$ is the generator such that $(X^i_t, X^i_t a) \in M_i$. The event $A_i$ occurs if (but not only if)
(1) there is no $M_i$-refresh time $t < T_i$ for $X^i$ and
(2) there is no $M_i$-refresh time $t > T_i$ for $X^i$.
Here, (1) guarantees that (3) $X^i_{T_i^+} \ne X^{i+1}_{T_i^+}$, while (2) guarantees that if (3) holds, then $X^i_{t} \ne X^{i+1}_{t}$ for all $t > T_i$. 
Note that the law of $(X^i_t)_{0 \le t \le T_i}$ is the same as the law of
\vadjust{\kern2pt}%
$\bigl((X^i_{T_i})^{-1} X^i_{T_i - t}\bigr)_{0 \le t \le T_i}$ by reversibility. In addition, the reflection in $M_i$ maps a path on one side of the wall to a path on the other side, preserving the times on the wall. Therefore, the probability of (1) is at least the probability of (2). The events (1) and (2) are also independent.
Hence, $\P(A_i) \ge \eta(W, S; \rt)^2$, which yields the desired inequality $\sigma(\rt') - \sigma(\rt) \ge 2 \epsilon \eta(W, S; \rt)^2 > 0$ in light of \rref l.norefresh/. 
\Qed

Combining Proposition \ref{prop:productrates} with Remark \ref{rem:finaffrates} and Theorem \ref{the:hyprates} gives Theorem \ref{the:strictgenerators}. 

\begin{Remark}
The assumption in Theorem \ref{the:strictgenerators} that each rate is positive is necessary --- if enough of the rates are zero, then the escape speed may be zero, even when one of the rates is increased. In a scenario where some of the rates are zero, one should just remove those generators, and apply Theorem \ref{the:strictgenerators} to the resulting smaller Coxeter system. 
Of course, if a rate is zero, increasing that rate will strictly increase the speed if the same holds when that rate is positive.
\end{Remark}

\begin{Remark}
The first part of Theorem \ref{the:hyprates} can be proved in the case of trees by using an implicit formula for the speed:
Namely, \cite{SawSte} gives a formula for the speed of the embedded discrete-time random walk. Multiplying that speed by the sum of the rates gives the following formula for the continuous-time speed. For a free product of $p \ge 3$ copies of $\Z/2\Z$ and corresponding rates $r_1, \dots, r_p > 0$, write $\soltn$ for the unique positive solution to the equation $\sum_{i=1}^p \bigl(\sqrt{\soltn^2 + r_i^2} - r_i\bigr) = (p-2) \soltn$. (Existence and uniqueness follow from the fact that the left-hand side is a convex function of $\soltn$ that passes through $(0, 0)$ with derivative 0 and is asymptotic to $p \soltn$ as $\soltn \to\infty$.)
The continuous-time speed $\sigma(\rt)$ is then $\soltn^{-1} \sum_{i=1}^p r_i \bigl(\sqrt{\soltn^2 + r_i^2} - r_i\bigr)$, which can also be written as $\sum_{i=1}^p \frac{r_i \soltn}{\sqrt{\soltn^2+r_i^2} + r_i}$. 
Because $r_j \mapsto \sqrt{x^2 + r_{\smash{j}\phantom{i}}^2} - r_j$ is strictly decreasing for fixed $x > 0$, it follows that $\partial \soltn/\partial r_j > 0$. Now
\[
\frac{\partial \sigma(\rt)}{\partial r_j}
=
\frac{\partial \soltn}{\partial r_j} \frac{1}{\soltn} \Bigl(- \sigma(\rt) + \sum_{i=1}^p \frac{r_i \soltn}{\sqrt{\soltn^2 + r_i^2}} \Bigr)
+
\frac{r_j}{\soltn} \Bigl({\textstyle \sqrt{(\soltn/r_j)^2  + 1}} - 2 + \frac1{\sqrt{(\soltn/r_j)^2  + 1}} \Bigr).
\]
From the second formula for $\sigma(\rt)$, it is clear that the first term in parentheses is strictly positive; it is elementary that the second term in parentheses is also strictly positive. Thus, $\sigma(\rt)$ is indeed strictly increasing in each rate.
\end{Remark}

\subsection{Further questions}

The proof of Proposition \ref{prop:extrars} gives rise to the following question.

\begin{Question}\label{q:couplefarther}
Given two random walks, with one having higher rates than the other, is it possible to couple them so that the more active walk is always at least as far from the initial position as the other?
\end{Question}

It might be plausible that a slightly cleverer construction in the proof of Proposition \ref{prop:extrars} could produce such a coupling. Surprisingly, this is not possible in even the simplest of cases.

\begin{Proposition} \label{prop:cantforever}
Question \ref{q:couplefarther} is impossible in the setting of Example \ref{ex:treecayley} even in cases where $r_b = r_c = 0$. In such cases, the walker is moving back and forth between just two states, at the rate $r_a$. Similarly, it is not always possible for $r_a, r_b, r_c > 0$.
\end{Proposition}
\begin{proof}
Consider two such walks, one with $r_a = 1$, and the other with $r_a = 100$, and the behaviour of these over the first unit of time. The first process has a significant chance of, for instance, staying at the initial vertex until a time between $0.4$ and $0.6$, moving to the other vertex, and staying there until time $1$. The more active process is very unlikely to take a path that remains farther from the starting vertex than such a path, because that would require staying at the other vertex between times $0.6$ and $1$, which is unlikely for such an active process. 

If $r_b, r_c > 0$, then one can choose a sufficiently small $t$ and make two choices of $r_a$, the smaller one being equal to $1/t$, so that a similar analysis holds.
\end{proof}

\procl p.recurrent-coupling
The coupling of Question \ref{q:couplefarther} is not possible for two walks on a recurrent graph, one of which has edges all with rate $\rho_1$, and the other all with rate $\rho_2 > \rho_1 > 0$.
\endprocl

\rproof
Let both random walks start at $o$. If such a coupling were possible, then we could similarly couple random walks with rates $(c \rho_1, c \rho_2)$ for any $c > 0$. In particular, we could couple with rates $\bigl(\rho_1 (\rho_2/\rho_1)^n, \rho_2 (\rho_2/\rho_1)^n\bigr)$ for every nonnegative integer $n$. Combining such couplings, we could couple with rates $(\rho_1, c \rho_2)$ for arbitrarily large $c$. By recurrence, the times between visits of the random walk with rate $\rho_2$ to the vertex $o$ have a distribution that is dominated by a geometric sum of exponential random variables with rate $\rho_2$, i.e., by an exponential random variable, which means that given $\epsilon > 0$, for all sufficiently large $c$ the probability is at least $1 - \epsilon$ that $o$ is visited in every interval of length $\epsilon$ in $[0, 1]$ by the random walk with rate $c \rho_2$. By the coupling, the same is true for the random walk with rate $\rho_1$, which forces that random walk to be at $o$ a.s.\ during the entire interval $[0, 1]$, a contradiction. 
\Qed

Question \ref{q:couplefarther} appears to be difficult to decide even in other simple settings.

\begin{Question}\label{q:onz}
On which transient Cayley graphs is the coupling of Question \ref{q:couplefarther} possible for two walks, one of which has edges all with rate $\rho_1$, and the other all with rate $\rho_2 > \rho_1 > 0$? 
\end{Question}

It is even possible to ask a version of this question for a chain with only three states.

\begin{Question}\label{q:on3}
Is the coupling of Question \ref{q:couplefarther} possible for two walks on the set $\{0,1,2\}$ that both start at $0$, where walk $i$ moves from $0$ to $1$, $1$ to $0$, and $1$ to $2$ with rate $\rho_i$? Note that once the walker reaches the state $2$, it stays there.
\end{Question}

Neither Question \ref{q:onz} nor \ref{q:on3} can be answered by a Markovian coupling, because it would be possible for both walks to be at a neighbor of the starting vertex, and then the more active walk has a greater chance to move back to the starting vertex than the less active walk does.

Finally, here are two more questions whose answers are unknown to us.

\procl q.concave
Is the escape speed concave in the rates for Coxeter systems?
\endprocl

\procl q.fields
Suppose that we have two random fields of rates on a Coxeter system with the first being at most the second a.s. Suppose that the law of each field is invariant under left group multiplication. Is the expected speed of the first at most that of the second?
\endprocl

\section{Approach to stationarity}
\label{sec:stationarity}

We will now consider random walks on finite groups and how the distance of the law of $Z_t$ to the stationary distribution changes when the random walk rates are increased.

\begin{Question}\label{q:closer}
Given a continuous-time random walk $(Z_t)_t$ on a finite group with positive rates on the generators and a time $t$, consider the distance between the distribution of $Z_t$ and the stationary distribution. Does the distance necessarily decrease when a random walk rate is increased? 
\end{Question}

It seems natural to think that the distance would decrease, as increasing a rate might be thought of as injecting additional randomness into the walk. In this section we will show that the answer is yes in several families of nice examples, and then give several examples where the answer is no in general, including some demonstrating rather unusual behaviour. Furthermore, there appear to be interesting patterns in the sorts of examples that we are able to produce.

There are several different measures of the distance from stationarity. Of course, the stationary distribution is uniform.

\begin{Definition}
For $p \in [1, \infty)$, the \defn{$\ell^p$-distance} between distributions $f$
and $g$ is 
\[
\ell^p(f,g) :=
\biggl(\sum_v|f(v)-g(v)|^{p}\biggr)^{\frac{1}{p}},
\]
where the sum is
taken over the states $v$ of the chain. The \defn{$\ell^\infty$-distance} between distributions $f$ and $g$ is $\ell^\infty(f,g) := \max_v |f(v) - g(v)|$.
The \defn{Hellinger distance} between distributions $f$ and $g$ is 
\[
\biggl(\frac12 \sum_v\bigl|\sqrt{f(v)}-\sqrt{g(v)}\bigr|^{2}\biggr)^{\frac{1}{2}}.
\]
\end{Definition}

\begin{Definition}
The \defn{entropy} of a distribution $f$ is the sum over states $v$ of $-f(v)\log f(v)$.
\end{Definition}

Observe that because the function $x\log x$ is convex on $(0, 1)$, the distribution of maximum entropy is the uniform distribution. Thus in this context, Question \ref{q:closer} asks whether increasing one of the rates causes the entropy of $Z_t$ to increase. 

One particularly strong way in which Question \ref{q:closer} may be answered in the affirmative would be if increasing one of the rates results in decreasing the vector of transition probabilities in the majorization order.

\procl d.major
Let $f$ and $g$ be vectors of length $n$. Denote by $f_{[1]} \ge f_{[2]} \ge \cdots \ge f_{[n]}$ the decreasing rearrangement of $f$, and similarly for $g_{[i]}$. We say that $f$ \defn{majorizes} $g$ if for each $i$, $$f_{[1]} + \cdots + f_{[i]} \geq g_{[1]} + \cdots + g_{[i]},$$ with equality when $i=n$ (which is automatic if $f$ and $g$ are probability vectors). Majorization defines a partial order on probability vectors, with the largest elements being the extreme points such as $(1,0,0,\dots,0)$ and the smallest element being $(\frac{1}{n},\frac{1}{n},\dots,\frac{1}{n})$. 
\endprocl

Recall the inequality of Hardy--Littlewood--P\'olya--Karamata that if $f$ majorizes $g$ and $\phi$ is convex, then $\sum_{i=1}^n \phi(f_i) \ge \sum_{i=1}^n \phi(g_i)$; if $\phi$ is strictly convex, then equality holds only if $f = g$. For example, if $f$ and $g$ are probability vectors and $\phi(x) = |x - 1/n|^p$ for some $p \ge 1$, then we may conclude that the $\ell^p$-distance between $f$ and the uniform distribution is at least the $\ell^p$-distance between $g$ and the uniform distribution. By letting $p \to\infty$, we get the same for the $\ell^\infty$-distance. A similar inequality holds for Hellinger distance and, with opposite sign, for entropy. Recall also that $f$ majorizes $g$ iff $g$ is gotten by applying a doubly stochastic matrix to $f$, such as by replacing $(f_i, f_j)$ by $\bigl(\alpha f_i + (1-\alpha)f_j, (1-\alpha)f_i + \alpha f_j\bigr)$, where $\alpha \in [0, 1]$. This is equivalent to moving some part of the larger of $f_i$ and $f_j$ to the smaller of the two. More generally, doubly stochastic matrices are convex combinations of permutation matrices.

In particular, for continuous-time random walk on any finite graph, the distribution of $Z_t$ decreases in the majorization partial order as $t$ increases.
In some cases, increasing one of the rates causes the distribution $Z_t$ to decrease in the majorization partial order, which implies that the $\ell^p$-distance to uniform has decreased for each $p \ge 1$ and that the entropy has increased. 

\subsection{Positive examples}

There are several special cases where the answer to Question \ref{q:closer} is positive. 

\subsubsection{Coxeter systems}

Throughout this subsection, we work in a Coxeter system, $(W, S)$. We will show that increasing the generator rates has the effect of decreasing in the majorization order the probability distribution of $Z_t$ at any fixed time $t$, which therefore implies decreasing distance to the uniform distribution in all senses mentioned. In order to prove this, we will use the \defn{Bruhat order}, the partial order $\le$ on $W$ defined by taking the transitive and reflexive closure of the relation 
\[
\bigl\{(v, Lv) \st L \textrm{ is a reflection and } |v| < |Lv|\bigr\}.
\]
We first extend Theorem \ref{the:coxmonotoneprobs} to the Bruhat order:

\procl t.bruhatmono
Let $x < y$ in the Bruhat order and $t > 0$. Then $p_t(o, x) > p_t(o, y)$.
\endprocl

\rproof 
It suffices to prove this when $y$ is a reflection $L$ of $x$ in a wall $M$. Let $T$ be the first $M$-refresh time of the random walk. For every path with $Z_t = y$, we have $T < t$ in light of Lemma \ref{lem:walldef}. If we change the coin flip at time $T$, then the walk will end instead at $L y = x$. This defines a probability-preserving bijection of the set of paths and coin flips (on $[0, t]$) from $o$ to $y$ to the set of paths and coin flips from $o$ to $x$ that have $T < t$. Since there is a positive-probability set of paths from $o$ to $x$ during which there is no $M$-refresh time, we get the desired strict inequality.
\Qed

We now look at how partial orders relate to majorization, based on \cite{BD}.

\procl d.poset
Let $\lep$ denote a partial order on a set, $A$. A function $f\colon A \to \R$ is called \defn{decreasing} if $f(x) \ge f(y)$ whenever $x \lep y$. A subset $B \subseteq A$ is called \defn{decreasing} if its indicator $\I B$ is decreasing. If $f$ and $g$ are two functions on $A$, then $f$ \defn{$\lep$-majorizes} $g$ if $f$ and $g$ are decreasing and $\sum_{x \in B} f(x) \ge \sum_{x \in B} g(x)$ for all decreasing sets $B$, with equality when $B = A$.
\endprocl

\procl l.posetmajor
Let $\lep$ denote a partial order on a finite set, $A$. If $f$ $\lep$-majorizes $g$, then $f$ majorizes $g$.
\endprocl

\rproof
For $i < |A|$, there is a decreasing set $B$ of cardinality $i$ on which $g$ attains its $i$ largest values. The sum of the $i$ largest values of $f$ is at least the sum of $f$ over $B$, which, by definition, is at least the sum of $g$ over $B$.
\Qed

\procl c.bruhatmajor
Let $\rt$ and $\rt'$ be two sets of rates on $S$ and $t$ and $t'$ be two positive times. Let $Z$ and $Z'$ be the corresponding random walks on $W$. Suppose there is a coupling of $Z_t$ and $Z'_{t'}$ such that $Z_t \le Z'_{t'}$ a.s.\ in the Bruhat order. Then the law of $Z_t$ majorizes the law of $Z'_{t'}$.
\endprocl

\rproof
Let $f$ be the law of $Z_t$ and $g$ be the law of $Z'_{t'}$. By \rref t.bruhatmono/, both $f$ and $g$ are decreasing. If $B$ is a decreasing subset of $W$, then $\P[Z_t \in B] \ge \P[Z'_{t'} \in B]$ because of the coupling, whence $f$ $\le$-majorizes $g$. Thus, the result follows from \rref l.posetmajor/. 
\Qed

\procl t.ratemajorCox
Let $\rt$ and $\rt'$ be two sets of rates on $S$ with $r_s \le r'_s$ for all $s \in S$. Let $t > 0$. Denote the corresponding transition probabilities by $p_t(x, y; \rt)$ and $p_t(x, y; \rt')$. If\/ $W$ is finite, then $p_t(o, \cdot; \rt)$ majorizes $p_t(o, \cdot; \rt')$ with inequality if $\rt \ne \rt'$.
\endprocl

\rproof
Suppose that $\rt = \rt'$ except for $r_s < r'_s$ for one specific $s$. Consider how many extra rings in $[0, t)$ there are of $R_s$ with the rates $\rt'$. If there is only one, then the argument of Corollary \ref{cor:coxmonotonespeed} showed, with the same notation, that we may couple $Z^1_t$ and $Z^2_t$ so that $Z^1_t \le Z^2_t$ in the Bruhat order, and inequality holds with positive probability. 
By induction, this same conclusion extends to any finite set of additional rings of $R_s$. Therefore, we may couple the two random walks so that $Z^1_t \le Z^2_t$ conditional on the extra rings of $R_s$. Hence we may couple them so that this inequality holds without such conditioning. By \rref c.bruhatmajor/, we deduce that the distribution of $Z^1_t$ majorizes that of $Z^2_t$. For the general case of rates, we may change the rates one by one until $\rt$ becomes $\rt'$, still yielding majorization. 
\Qed

Of course, there is a superficial generalization of this result as in \rref r.semicox/. With the proper definition, this result also holds for infinite $W$.
Similarly to \rref t.coxdiscretedist/, we may deduce from our arguments the following, which also easily implies \rref t.ratemajorCox/:

\procl t.coxdiscretemajor
Let $(W, S)$ be a Coxeter system. Let $\seqs$ and $\seqs'$ be finite sequences from $S$ of lengths $n$ and $n'$, respectively, with $\seqs$ a proper subsequence of $\seqs'$. 
Let $\cfs$ be a Bernoulli$(1/2)$ process.
Then the probability distribution of $\cfm(\seqs, \cfs_{n})$ strictly majorizes that of $\cfm(\seqs', \cfs_{n'})$. 
\qed
\endprocl

\procl r.stoch-domin
There are additional consequences of the argument used in \rref t.bruhatmono/.
Suppose that $M$ is a wall. Write $M^+$ for the set of vertices $v$ for which there is a path joining $v$ to $o$ without using an edge of $M$, and write $M^-$ for the remainder. Since the Cayley graph of $(W, S)$ is connected, Lemmas \ref {Lemma:crossing} and \ref{lem:walldef} imply that the reflection $L_M$ in $M$ interchanges $M^+$ and $M^-$. Let $\tau_{M^-} := \inf \{t > 0 \st Z_t \in M^-\}$. A direct analogue of the reflection principle (for one-dimensional random walks and Brownian motion) is that for all $t > 0$ and all $A \subseteq M^+$,
\[
\P[\tau_{M^-} < t,\ Z_t \in A]
=
\P[Z_t \in L_M A].
\]
Therefore,
\[
\P[\tau_{M^-} < t]
=
2 \P[Z_t \in M^-].
\]
Write $\tau_v := \inf \{t > 0 \st Z_t = v\}$. If $x < y$, then $\tau_y$ strictly stochastically dominates $\tau_x$, that is, for all $t > 0$, $\>\P[\tau_x < t] > \P[\tau_y < t]$. In particular, $\P[\tau_x < \infty] \ge \P[\tau_y < \infty]$. Let $\lambda_v(t)$ denote the Lebesgue measure of the set of times in $[0, t]$ when the random walk is at $v$. As a consequence of the strict stochastic domination inequality, we obtain that $\lambda_x(t)$ strictly stochastically dominates $\lambda_y(t)$ when $x < y$ and $t > 0$.
We leave the arguments to the reader.
\endprocl

\subsubsection{Abelian groups and conjugacy classes}

In an abelian group, the answer to Question \ref{q:closer} is always positive, because increasing a rate results in extra multiplications by random group elements, and in an abelian group, we may consider these extra multiplications to take place at the end of the random walk. 

As a generalisation of this observation, consider a random walk on a group $G$ with a generator $g$ every conjugate of which is also a generator. If all of the rates on this conjugacy class are increased the same amount, then this results in multiplying by additional group elements partway through the walk. Those extra elements are uniformly distributed in the conjugacy class of $g$, independently of all other steps. This is equivalent to multiplying by the same number of independent, uniformly chosen conjugates of $g$ at the end.  Thus, this also moves the resulting distribution down in the majorization order.

\subsubsection{All groups with special distances}

As long as we measure distance to stationarity in special ways, increasing rates will always decrease distance on any group.

\begin{Proposition}
\label{prop:pinf}
For any finite group, increasing any of the rates always decreases the $\ell^\infty$-distance and the $\ell^2$-distance to the uniform distribution.
\end{Proposition}
\begin{proof}
This relies on some well-known calculations. We have
\begin{align} \label {eq:infty}
\begin{split}
|p_t(x, y) - 1/n|
&=
\Bigl|\sum_z \bigl(p_{t/2}(x, z) - 1/n\bigr)\bigl(p_{t/2}(z, y) - 1/n\bigr)\Bigr|
\\ &=
\Bigl|\sum_z \bigl(p_{t/2}(x, z) - 1/n\bigr)\bigl(p_{t/2}(y, z) - 1/n\bigr)\Bigr|
\\ &\le
\sqrt{\sum_z \bigl(p_{t/2}(x, z) - 1/n\bigr)^2\sum_z \bigl(p_{t/2}(y, z) - 1/n\bigr)^2}
\\ &=
\sqrt{\bigl(p_{t}(x, x) - 1/n\bigr)\bigl(p_{t}(y, y) - 1/n\bigr)}
\\ &=
p_{t}(o, o) - 1/n.
\end{split}
\end{align}
Note that $p_t(o, o) \ge 1/n$ because the first two lines above, without absolute values, show that 
\begin{equation} \label {eq:2infty}
p_t(o, o) - 1/n = \sum_z \bigl(p_{t/2}(o, z) - 1/n\bigr)^2 \ge 0.
\end{equation}
In particular, \eqref{eq:infty} gives
\[
\max_{x, y} |p_t(x, y) - 1/n|
=
p_{t}(o, o) - 1/n.
\]
When any rate is increased, $p_t(o, o)$ decreases, as we noted in the introduction. Thus, the $\ell^\infty$-distance decreases. In light of \eqref{eq:2infty}, the $\ell^2$-distance at time $t$ equals the square root of the $\ell^\infty$-distance at time $2t$, whence it is also decreasing in the rates.
\end{proof}

\subsection{Negative examples} \label{sec:negative}

Symmetric groups give Coxeter systems when generated by adjacent transpositions. Changing generators, however, may yield entirely different behaviour. Indeed,
every finite group is a subgroup of a symmetric group. Therefore, any counterexample can be exhibited on a symmetric group by adding generators to those for the subgroup and making their rates 0 or very close to 0. However, symmetric groups being very large means that it can be hard to find examples by searching in symmetric groups.

In this section, many of our examples were found by random numerical searching and numerical calculation --- choosing sets of random generators and rates and exploring the consequences of increasing one of those rates. As we move through the different types of examples, we will attempt to give some idea of how difficult it was to find them, as an indication of the frequency of similar examples. We will usually be considering groups with up to about $120$ elements, and numbers of random samples in the hundreds of thousands. In cases where we were unable to find examples, we will indicate when we spent enough effort that we were surprised not to find them.

In most cases, we used sets of four generators, with rates between about $1$ and $10$. Indeed, when choosing five generators at random, the most extreme examples often either included the identity or two copies of the same generator, suggesting that it is easier to find these examples with four distinct generators than with five.

For many of the groups we examined, it was not difficult to find examples where increasing a rate increased the $\ell^1$-distance from stationarity. This includes the symmetric groups $S_4$ and $S_5$, dihedral groups with as few as ten elements, and dicyclic groups with as few as twenty elements. We did not find these examples for split metacyclic groups. We will restrict our discussion now to dihedral and symmetric groups, as two quite different families of groups.

For the dihedral group $D_{11}$ with three involutions as generators, we found an example where the entropy was not monotone in the rates. Note that dihedral groups are Coxeter groups, with Coxeter generators being two involutions.

As we increase $p$ from $1$ towards $2$, it becomes more difficult to find examples where the $\ell^p$-distance increases, which is consistent with Proposition \ref{prop:pinf}. For instance, in the dihedral group $D_5$, we can find examples for values of $p$ between $1$ and $1.4$, and in the group $D_7$, for values of $p$ between $1$ and $1.8$. It also seems easier to find these kinds of examples in dihedral groups $D_n$ when $n$ is prime compared to a composite number of similar size, such as $D_{31}$ or $D_{41}$ compared to $D_{30}$ or $D_{40}$.

While Proposition \ref{prop:pinf} says that there are no examples of this kind for $p=2$, it is possible to get very close to this value. For instance, it was not uncommon for us to find examples exhibiting this behaviour for $p$ up to $1.99$ or from $2.01$, and individually optimised examples worked for $p$ as large as $1.9999$.

The next natural question is what happens for values of $p$ larger than two. In the case of dihedral groups, we found examples of sets of generators and rates which increase the $\ell^p$-distance from uniform for a range of values of $p$ between $2+\eps$ and $4-\eps$, for $\eps$ as small as $0.001$, but we were unable to find any examples where the $\ell^p$-distance increases for $p=4$, searching in dihedral groups as large as $D_{131}$. This seems rather surprising. A generic example of this sort has the $\ell^p$-distance increasing for $p$ in some interval contained in $(2, 4)$. When the rates are adjusted so that the right endpoint of this interval moves closer to $p=4$, the left endpoint seems to move closer to $p=2$, and vice versa.

Likewise, we found an example in the dihedral group $D_{41}$ for which the $\ell^p$-distance increases for values of $p$ between $4.001$ and $5.995$. Interestingly, this example also increases the $\ell^p$-distance for $p$ up to $p=1.997$.

We were unable to find dihedral examples either at $p=6$ or for larger values of $p$. Proposition \ref{prop:pinf} explains why distance cannot increase for $p=2$, but we have no similar explanation for $p=4$ or $p=6$. While the $\ell^p$-distance may increase for ranges of $p$ less than $2$ or between $2$ and $4$, we were unable to find any dihedral examples exhibiting this behaviour in both ranges simultaneously.

The proof of Proposition \ref{prop:pinf} relates the behaviour of the $\ell^2$-distance at time $t$ to the $\ell^\infty$-distance at time $2t$. One might hope, then, that if there are examples where the $\ell^p$-distance increases for $p$ near $2$, then there should be similar examples where the $\ell^p$-distance increases for $p$ near $\infty$ --- that is, for $p$ very large --- at twice the time. Yet we were unable to find such dihedral examples.

We looked for similar examples to these in the symmetric group $S_5$. As with the dihedral groups, we found examples where the $\ell^p$-distance increased for all $p$ up to a value quite close to $2$, and examples where it increased for most $p$ between $2$ and $4$, but also new types of behaviour, including one where the $\ell^p$-distance increased for $p$ between $2.004$ and $4.02$, and several where it increased for $p$ from about $3$ up to as large as $300$.

We did not find examples where the distance increased at $p=4$ in dicyclic groups.

\section{A ray}

We now turn our attention from Cayley graphs to 
an infinite one-ended path, and allow arbitrary rates on the edges. We label vertices by the nonnegative integers, and each pair of consecutive numbers is connected by an edge. For each positive integer $i$, let $r_i$ be the (nonnegative) rate of the Poisson clock on the edge between $i-1$ and $i$. Let $i_0 := \inf \{ j \st r_j = 0\} \in [1, \infty]$. For simplicity, we assume that the rates are such that explosions do not occur. We will see several interesting phenomena in this case. The random walk will begin from $0$.

\begin{Proposition}\label{prop:decreasingprob}
For each time $t > 0$, the probability $p_t(0, i)$ that the walker is at location $i$ is a strictly decreasing function of $i$ for those $i < i_0$.
\end{Proposition}
\begin{proof}
We first prove that $p_t(0, i)$ is weakly decreasing in $i$.

Let $f_t(i, j)$ be the probability that a random walk starting from $i$ visits $j$ at some time before $t$.
Because there are no explosions, $\lim_{N \to\infty} f_t(0, N) = 0$. Consider the sequence of chains whose rates are $r_i$ for $i \le N$ and $0$ for $i > N$; write $p_t^N(i, j)$ for their transition probabilities. 
Since $|p_t(0, i) - p_t^N(0, i)| \le f_t(0, N)$, it follows that $\lim_{N \to\infty} p_t^N(0, i) = p_t(0, i)$ for all $i$. Therefore, it suffices to prove the claim for a chain where only finitely many rates are nonzero. We now assume that condition.

Consider the refresh times $R_i$ of Definition \ref{def:refresh}. There are only finitely many such times before $t$ a.s. If we condition on the sequence of all refresh rings that occur between time $0$ and time $t$, ordered by time, then the probability distribution function of the walker's location at time $t$ is a nonincreasing function of $i$. This is because the coin flip, or randomization, at each time of $R_i$ acts on the walker's distribution function at that time by averaging the probabilities at $i-1$ and $i$, and this operation preserves monotonicity. The walker is initially at position $0$ with probability $1$, and this initial distribution function is monotone.

Integrating these conditional distribution functions with respect to the distribution of all refresh rings before time $t$ completes the proof of the claim. 

Next we show that the transition probabilities are strictly decreasing. We no longer assume that only finitely many rates are nonzero.
Let $\tilde p_t(0, \cdot)$ denote the probability distribution conditional on the set of refresh rings before time $t$.
There is a collection of positive probability of sets of refresh rings for each of which $\tilde p_t(0, i-1) > \tilde p_t(0, i)$ for $i < i_0$. Because $\tilde p_t(0, i-1) \ge \tilde p_t(0, i)$ for each possible set of refresh rings and $p_t(0, i)$ is the expectation of $\tilde p_t(0, i)$, it follows that $p_t(0, i)$ is strictly decreasing for $0 \le i < i_0$.
\end{proof}

This analysis also lets us discuss the speed at which the walker moves away from $0$.

\begin{Lemma}\label{lem:onemorer}
In the language of the proof of Proposition \ref{prop:decreasingprob}, consider the following two distribution functions on $\Z^+$. The function $f_1$ is defined by starting with a unit mass at $0$ and applying any fixed finite sequence of refresh rings. The function $f_2$ is defined similarly, using the same sequence of refresh rings, except with one additional refresh ring $R_j$ occurring partway through the sequence. Then $f_2$ dominates $f_1$, that is, $\sum_{i = 0}^k f_1(i) \ge \sum_{i = 0}^k f_2(i)$ for every $k \in \Z^+$.
\end{Lemma}
\begin{proof}
Consider the evolution of the difference $f_2 - f_1$. This is zero until the extra $R_j$ is applied to $f_2$. At this point, $f_2 - f_1$ is zero at all points except $j-1$ and $j$. At $j-1$ it is negative (or zero), at $j$ it is positive (or zero), and these two differences have the same size (but opposite sign). 

At this point, the difference $f_2 - f_1$ has the property that when summed over the states $0$ to any $k$, it is nonpositive. This property is preserved by the application of any of the $R_i$, which completes the proof.
\end{proof}

This implies that the location of the random walk at each time is stochastically strictly increasing in each of the rates $r_j$:

\begin{Corollary}\label{cor:cumulprobs}
At every time $t$ and nonnegative integer $i < i_0$, the probability that the walker is at a position between $0$ and $i$ is a strictly decreasing function of each $r_j$ for $j < i_0$.  
\end{Corollary}
\begin{proof}
As in the proof of Proposition \ref{prop:decreasingprob}, we first prove the statement without the strictness under the assumption that $i_0 < \infty$.
When some $R_j$ rings one additional time, then the result follows from Lemma \ref{lem:onemorer}. Applying this fact repeatedly gives the result. It then follows even if $i_0 = \infty$, but without the strictness part. 

To show strictness, it again suffices, for each $i$, to exhibit a collection of positive probability of sequences of refresh rings with a marked additional ring of $R_j$ that each give a strict inequality. Indeed, such a collection is formed by the sequences that when restricted to the refresh rings for edges $k \le i \vee j$ are sequences of the form $(R_1, R_2, \ldots, R_{i \wedge j}, \ldots, R_{i \vee j})$, where (one of) the $R_j$ is marked additional.
\end{proof}

\begin{Corollary}\label{cor:expdist}
At every time $t$, the expected distance of the walker from $0$ is a strictly increasing function of each $r_j$ for $j < i_0$.
\end{Corollary}
\begin{proof}
The expected distance from $0$ is the sum over all $i$ of the probabilities that the walker is farther away than $i$. Hence the result follows from Corollary \ref{cor:cumulprobs}.
\end{proof}

We are also interested in how much time the walker spends at each vertex.

\begin{Proposition}\label{prop:timedom}
For our random walk on the infinite one-ended path, let $t$ be an arbitrary time. Then the time spent at $0$ between time $0$ and time $t$ stochastically dominates the time spent at $1$.
\end{Proposition}
\begin{proof}
As in previous proofs, we may assume that $i_0 < \infty$.
Without loss of generality, rescale time so that the rate $r_1$ is equal to $1$. Fix $t > 0$. We define recursively the following sequences of times. Let $A_1 := 0$. Then let $B_k$ be the infimum of the times in $[0, t]$ after time $A_k$ at which the walker is at the vertex $1$ if there is such a time and $B_k := t$ if not, and $A_{k+1}$ be the infimum of the times in $[0, t]$ after $B_k$ at which the walker is at $0$ if there is such a time and $t$ if not, for each $k$. For each $k$, let $X_k$ be the amount of time the walker spends at vertex $0$ between time $A_k$ and $A_{k+1}$, and $Y_k$ be the amount of time spent at $1$ between $B_k$ and $B_{k+1}$. 

If we didn't stop counting at time $t$, then both $X_k$ and $Y_k$ would be distributed as exponential random variables with rate $1$, because they count the time spent at vertex $0$ or $1$, respectively, until the Poisson clock corresponding to the edge $(0, 1)$ rings while the walk is at one of its endpoints, and this clock has rate $1$. Taking this into account, we see that the conditional distribution of $X_k$ given $A_k$ and $(X_i, Y_i)$ for all $i < k$ is the same as the distribution of $\min\bigl(\Exp(1),t-A_k\bigr)$. We can't give such a concrete expression for the conditional distribution of $Y_k$, because the walker may move to $2$ and beyond before moving back to $1$ and then to $0$. We can say, though, that $Y_k$ given $A_k$, $X_k$, and $(X_i, Y_i)$ for all $i < k$ is stochastically dominated by $\min\bigl(\Exp(1),t-A_k-X_k\bigr)$. In particular, this means that $X_k$ dominates $Y_k$ given $(X_i, Y_i)$ for all $i < k$. 
It follows that the sum of the $X_k$ stochastically dominates the sum of the $Y_k$.
Since the total time spent at $0$ between $0$ and $t$ is the sum of the $X_k$, and the time spent at $1$ is the sum of the $Y_k$, we obtain the desired result.
\end{proof}

In the preceding proof, note that we gave an explicit expression for $X_k$, and an upper bound for $Y_k$. If the walker was allowed to move away from the vertex $0$ without moving to $1$, then we would only have obtained an upper bound for $X_k$, which would not allow us to compare $X_k$ and $Y_k$. Example \ref{ex:1to2} illustrates what may go wrong in such a scenario.

Surprisingly, even the simplest generalisations of Proposition \ref{prop:timedom} are not true. The following example shows that in the same setting, the walker need not spend more time at vertex $1$ than at $2$.

\begin{Example}\label{ex:1to2}
Fix $t > 0$.
Set the rate $r_1$ to be extremely large, so that the total time spent at vertex $0$ or $1$ before time $t$ is very likely to be almost the same. Take $r_2$ to be $1$ and $r_3$ to be zero (thus the walker will never reach any state beyond $2$, and the other rates are irrelevant). Then at time $t$, the probability that the walker has spent more than time $\frac{2t}{3}$ at vertex $1$ is almost zero, because it spends almost equal amounts of time at $0$ and $1$. However, there is probability bounded away from $0$ that the walker spends more time than $\frac{2t}{3}$ at state $2$ --- for instance, it could move to state $2$ by time $\frac{t}{4}$ and then the edge between states $1$ and $2$ might not ring again before time $t$.

Therefore, the time spent at $1$ does not stochastically dominate the time spent at $2$. 
\end{Example}


The next example illustrates how Corollary \ref{cor:expdist} may fail in cases that are only slightly more complicated. The state space will still be a path, but the walker will not start at an end.

\begin{Example} \label{ex:Zdist}
Expand the setting so that there is a vertex for each integer $i$, with $r_i$ still the rate for moving between $i-1$ and $i$. The walker still begins at $0$, but there are now possible states on either side of $0$.

If $r_i = 0$ for $i \ne 0, 1$ and $r_1 = 1$, then one can check by explicit calculation that $p_t(0, 0)$ is larger when $r_0 = 3$ than it is when $r_0 = 2$, whence the expected distance from $0$, being $1 - p_t(0, 0)$, is smaller when $r_0 = 3$ than it is when $r_0 = 2$.

For a more extreme example,
choose arbitrary positive integers $k$ and $n$, and set the rates on edges between $-k$ and $0$ to each be $1$, and rates on edges between $0$ and $n-1$ to be very large and the same. Rates outside these ranges are zero. At times that are much less than $1$, but large compared to the reciprocal of the large rate, the walker is approximately equally likely to be at any state between $0$ and $n-1$, for an expected distance of $\frac{n-1}{2}$. At large times, the walker is approximately equally likely to be at any of the $k+n$ states, for an expected distance of $\frac{n(n-1)+k(k+1)}{2(n+k)}$. For some values of $k$ and $n$, the expected distance from the starting state has decreased between these two regimes. 

If we take $k = \alpha n$ and consider the limit as $k$ and $n$ increase, the most extreme ratio occurs when $\alpha = \sqrt{2} - 1$, which gives a ratio between the small-time and large-time expected distances of $\frac{1 + \sqrt{2}}{2}$.
\end{Example} 

One permissible generalisation of Proposition \ref{prop:timedom} is to settings where the excursions from the two vertices in question have the same distribution. This includes the setting of Cayley graphs with rates depending only on the generators.

\begin{Proposition}
Consider random walk on a graph whose edge rates are preserved by some group of graph automorphisms that acts transitively on the vertices. Then the time spent between $0$ and $t$ at the initial vertex $u$ stochastically dominates the time spent at each other vertex $v$.
\end{Proposition}
\begin{proof}
Let $\lambda(x,t)$ be the amount of time between $0$ and $t$ spent at the vertex $x$ by a random walk started at $x$ at time $0$. Because the graph is vertex-transitive, the law of this quantity does not depend on $x$. The duration $\lambda(x,t)$ is also increasing in $t$. Let $\tau$ be a random variable independent of the random walk and whose distribution is the same as that of the hitting time of the vertex $v$ when the random walk is started from $u$. Then we are comparing $\lambda(u,t)$ with $\lambda\bigl(v,(t-\tau)^+\bigr)$, whose law is the same as that of $\lambda\bigl(u,(t-\tau)^+\bigr)$. The random variable $\tau$ is positive, which completes the proof.
\end{proof}

\medskip
{\bf Acknowledgements: } We thank Kavita Ramanan for early discussions of some of these questions. We also thank Yiping Hu for asking us about convergence to stationarity. Thanks to Yuval Peres, who alerted us to \cite{PeresWinkler} and \cite{FillKahn} after seeing an earlier draft of our results. We also thank one of the referees for informing us of \rref b.HerKoz/. We are grateful to Matthias Weber for permission to include Figure \ref{fig:237cayley}.

\bibliographystyle{plain}
\bibliography{bib}

\end{document}